\documentclass[journal]{IEEEtran}
\usepackage{cite}
\usepackage{amsmath,amssymb,amsfonts, amsthm}
\usepackage{algorithmic}
\usepackage{graphicx}
\usepackage{textcomp}

\newtheorem{theorem}{Theorem}

\newtheorem{problem}{Problem}

\DeclareMathOperator*{\expect}{\mathbb{E}}
\title{Competitive Control}
\author{Gautam Goel and Babak Hassibi
\thanks{Gautam Goel is with the Department of Computing and Mathematical Sciences at Caltech (e-mail: ggoel@caltech.edu). }
\thanks{Babak Hassibi is with the Department of Electrical Engineering at Caltech (e-mail: bhassibi@caltech.edu).}}
\date{}

\begin{document}
\maketitle

\begin{abstract}
We consider control from  the  perspective  of competitive analysis.  Unlike much prior work on learning-based control, which focuses on minimizing regret against the best controller selected in hindsight from some specific class, we focus on designing an online controller which competes against a clairvoyant offline optimal controller.  A natural performance metric in this setting is competitive ratio, which is the ratio between the cost incurred by the online controller and the cost incurred by the offline optimal controller. Using operator-theoretic techniques  from  robust  control, we derive a computationally efficient state-space description of the the controller with optimal competitive ratio in both finite-horizon and infinite-horizon settings. We extend competitive control to nonlinear systems using Model Predictive Control (MPC) and present numerical experiments which show that our competitive controller can significantly outperform standard $H_2$ and $H_{\infty}$ controllers in the MPC setting.
\end{abstract}

\section{Introduction}
The central question in control theory is how to regulate the behavior of an evolving system which is perturbed by an external disturbance by dynamically adjusting a control signal. Traditionally, controllers have been designed to optimize performance under the assumption that the disturbance is drawn from some specific class of disturbances. For example, in $H_2$ control the disturbance is assumed to be generated by a stochastic process and the controller is designed to minimize the expected cost, while in $H_{\infty}$ control the disturbance is assumed to be generated adversarially and the controller is designed to minimize the worst-case cost. This approach suffers from an obvious drawback: if the controller encounters a disturbance which falls outside of the class the controller was to designed to handle, the controller's performance may be poor. In fact, the loss in performance can be arbitrarily large, as shown in \cite{doyle1978guaranteed}.

This observation naturally motivates the design of adaptive controllers which dynamically adjust their control strategy as they sequentially observe the disturbances instead of blindly following a prescribed strategy. The design of such controllers has attracted much recent attention in the online learning community (e.g. \cite{agarwal2019online, hazan2020nonstochastic, foster2020logarithmic}), mostly from the perspective of \textit{policy regret}. In this framework, the online controller is designed to minimize regret against the best controller selected in hindsight from some time-invariant comparator class, such as the class of static state-feedback policies or the class of disturbance-action policies introduced in \cite{agarwal2019online}. The resulting controllers are adaptive in the sense that they seek to minimize cost without making \textit{a priori} assumptions about how the disturbances are generated.

In this paper, we take a somewhat different approach to adaptive control: we focus on designing a controller which minimizes  the \textit{competitive ratio}
$$\sup_{w} \frac{ALG(w)}{OPT(w)},$$ where $ALG(w)$ is the control cost incurred by the online controller in response to the disturbance $w$ and $OPT(w)$ is the cost incurred by a \textit{clairvoyant offline optimal controller}. The clairvoyant offline optimal controller is the controller which selects the globally optimal sequence of control actions given perfect knowledge of the disturbance $w$ in advance; the cost incurred by the offline optimal controller is a lower bound on the cost incurred by any controller, causal or noncausal. A controller whose competitive ratio is bounded above by $C$ offers the following guarantee: the cost it incurs is always at most a factor of $C$ higher than the cost that could have been counterfactually incurred by any other controller, \textit{irrespective of the disturbance is generated}. Competitive ratio is a multiplicative analog of \textit{dynamic regret}; the problem of obtaining controllers with optimal dynamic regret was recently considered in \cite{goel2021regret, sabag2021regret, goel2021measurement}.

We emphasize the key distinction between policy regret and competitive ratio: policy regret compares the performance of the online controller to the best fixed controller selected in hindsight from some class, whereas competitive ratio compares  the performance of the online controller to the optimal dynamic sequence of control actions, without reference to any specific class of controllers. We believe the competitive ratio formulation of online control we consider in this paper compares favorably to the policy regret formulation in two ways. First, it is more general: instead of imposing \textit{a priori} some parametric structure on the controller we learn (e.g. state feedback policies, disturbance action policies, etc), which may or may not be appropriate for the given control task, we compete with the globally optimal clairvoyant controller, with no artificial constraints. Secondly, and more importantly, the controllers we obtain are more robust to changes in the environment. Consider, for example, a scenario in which the disturbances are picked from a probability distribution whose mean varies over time. When the mean is near zero, an $H_2$ controller will perform well, since $H_2$ controllers are tuned for zero-mean stochastic noise. Conversely, when the mean is far from zero, an $H_{\infty}$ controller will perform well, since $H_{\infty}$ controllers are designed to be robust to large disturbances. No fixed controller will perform well over the entire time horizon, and hence any online algorithm which tries to converge to a single, time-invariant controller will incur high cumulative cost. A controller which competes against the optimal dynamic sequence of control actions, however, is not constrained to converge to any fixed controller, and hence can potentially outperform standard regret-minimizing control algorithms when the environment is non-stationary.

%Our competitive ratio approach to control is partially inspired by series of works in online learning (e.g. \cite{herbster1998tracking, hazan2009efficient, jadbabaie2015online,  goel2019beyond}), which study online learning in non-stationary environments. In this setting, the goal is to obtain algorithms which compete with a dynamic sequence of actions instead of the best fixed action selected in hindsight; this formulation is natural when the reward-generating process encountered by the online algorithm varies over time. 

\subsection{Contributions of this paper}
We derive the controller with optimal competitive ratio, resolving an open problem in the learning and control literature first posed in \cite{goel2019online}. Our competitive controller is a drop-in replacement for standard $H_2$ and $H_{\infty}$ controllers and can be used anywhere these controllers are used; it also uses the same computational resources as the $H_{\infty}$-optimal controller, up to a constant factor. The key idea in our derivation is to reduce competitive control to $H_{\infty}$ control. Given an $n$-dimensional linear dynamical system driven by a disturbance $w$, we show how to construct a synthetic $2n$-dimensional linear system and a synthetic disturbance $w'$ such that the $H_{\infty}$-optimal controller in the synthetic system driven by $w'$ selects the control actions which minimize competitive ratio in the original system. 

We synthesize the competitive controller in a linearized Boeing 747 flight control system; in this system, our competitive controller obtains the competitive ratio 1.77. In other words, \textit{ it is guaranteed to incur at most 77\% more cost than the clairvoyant offline optimal controller, irrespective of how the input disturbance is generated}. Numerical experiments show that the competitive controller exhibits ``best-of-both-worlds" behavior, often beating standard $H_2$ and $H_{\infty}$ controllers on best-case and average-case input disturbances while maintaining a bounded loss in performance even in the worst-case. We also extend our competitive control framework to nonlinear systems using Model Predictive Control (MPC). Experiments in a nonlinear system show that the competitive controller consistently outperforms standard $H_2$ and $H_{\infty}$ controllers across a wide variety of input disturbances, often by a large margin.

Our results can be viewed as injecting adaptivity and learning into traditional robust control; instead of designing controllers which blindly minimize worst-case cost irrespective of the disturbance sequence they encounter, we show how to extend $H_{\infty}$ control to obtain controllers which dynamically adapt to the disturbance sequence by minimizing competitive ratio. 
\subsection{Related work}
Integrating ideas from machine learning into control has has attracted much recent attention across several distinct settings. In the ``non-stochastic control" setting proposed in \cite{hazan2020nonstochastic},  the online controller seeks to minimize regret against the class of disturbance-action policies in the face of adversarially generated disturbances. An $O(T^{2/3})$ regret bound was given in \cite{hazan2020nonstochastic}; this was improved to $O(T^{1/2})$ in \cite{agarwal2019online} and  $O(\log{T})$ in \cite{foster2020logarithmic}. These works 
focus on minimizing regret against a fixed controller from some parametric class of control policies (policy regret); a parallel line of work studies the problem of designing an online controller which minimizes regret against a time-varying comparator class (dynamic regret). Dynamic regret is a very similar metric to competitive ratio, which we consider in this paper, except that it is the difference between the cost of the online and offline controllers, rather than the ratio of the costs. The problem of designing controllers with optimal dynamic regret was studied in the finite-horizon, time-varying setting in \cite{goel2021regret}, in the infinite-horizon LTI setting in \cite{sabag2021regret}, and in the measurement-feedback setting in \cite{goel2021measurement}. Gradient-based algorithms with low dynamic regret against the class of disturbance-action policies were obtained in \cite{gradu2020adaptive,  zhao2021non}.

In this paper, we design controllers through the lens of \textit{competitive analysis}, e.g. we seek to design online algorithms which compete against a clairvoyant offline algorithm. This idea has a rich history in theoretical computer science and we refer to \cite{borodin2005online} for an overview. In \cite{goel2019online}, Goel and Wierman showed that competitive ratio guarantees in a narrow class of linear-quadratic (LQ) systems could be obtained using the Online Balanced Descent (OBD) framework proposed in \cite{chen2018smoothed}. A series of papers \cite{goel2019beyond, shi2020online} extended this reduction; a similar reduction was explored in \cite{goel2017thinking} in the context of multi-timescale control.  We emphasize that all prior work failed to obtain a controller with optimal competitive ratio, and relied on making nonstandard structural assumptions about the dynamics; for example, \cite{shi2020online} assumes that the disturbance affects the control input rather than the state. This paper is the first to obtain controllers with optimal competitive ratio in general LQ systems, in both finite-horizon and infinite-horizon settings.

\section{Preliminaries} \label{preliminaries-sec}

In this paper we consider the design of competitive controllers in the context of linear-quadratic (LQ) control. This problem is generally studied in two distinct settings: finite-horizon control in time-varying systems and infinite-horizon control in linear time-invariant (LTI) systems. We briefly review each in turn:

\textit{Finite-horizon Control}.
In this setting, the dynamics are given by the linear evolution equation
\begin{equation} \label{evolution-eq}
x_{t+1} = A_tx_t + B_{u, t} u_t + B_{w, t} w_t.
\end{equation}
Here $x_t \in \mathbb{R}^n$ is a state variable we seek to regulate, $u_t \in \mathbb{R}^m$ is a control variable which we can dynamically adjust to influence the evolution of the system, and $w_t \in \mathbb{R}^p$ is an external disturbance. We focus on control over a finite horizon $t = 0, \ldots, T-1$ and often use the notation $w = (w_0, \ldots, w_{T-1})$, $u = (u_0, \ldots, u_{T-1})$, $x = (x_0, \ldots, x_{T-1})$. We assume for notational convenience the initial condition $x_0 = 0$, though it is trivial to extend our results to arbitrary initialization. We formulate control as an online optimization problem, where the goal is to select the control actions so as to minimize the quadratic cost 
\begin{equation} \label{fh-lqr-cost}
\sum_{t=0}^{T-1} \left( x_t^*Q_t x_t + u_t^*R_t u_t \right),
\end{equation}
where $Q_t \succeq 0, R_t \succ 0$ for $t = 0, \ldots T - 1$.  We assume that the dynamics $\{A_t, B_{u, t}, B_{w, t}\}_{t=0}^{T-1}$ and costs $\{ Q_t, R_t\}_{t=0}^{T-1}$ are known, so the only uncertainty in the evolution of the system comes from the external disturbance $w$. For notational convenience, we assume that the system is parameterized such that $R_t = I$ for $t = 0, \ldots T-1$; we emphasize that this imposes no real restriction, since for all $R_t \succ 0$ we can always rescale $u_t$ so that $R_t = I$. More precisely, we can define $B'_{u, t} = B_{u, t}R^{-1/2}_t$ and $u_t' = R^{1/2}u_t$; with this reparameterization, the evolution equation (\ref{evolution-eq}) becomes $$x_{t+1} = A_tx_t + B_{u, t}' u_t' + B_{w, t} w_t,$$ while the state costs $\{Q_t\}_{t=0}^{T-1}$ appearing in (\ref{fh-lqr-cost}) remain unchanged and the control costs $\{R_t\}_{t=0}^{T-1}$ are all equal to the identity. This choice of parametrization greatly simplifies notation and is common in the control literature, see e.g. \cite{hassibi1999indefinite}.

\textit{Infinite-horizon Control}.
In this setting, the dynamics are given by the time-invariant linear evolution equation $$x_{t+1} = Ax_t + B_u u_t + B_w w_t,$$ where  $x_t \in \mathbb{R}^n$, $u_t \in \mathbb{R}^m$, and $w_t \in \mathbb{R}^p$. We focus on control over a doubly-infinite horizon and often use the notation $w = (\ldots, w_{-1}, w_0, w_1, \ldots)$, $u = (\ldots, u_{-1}, u_0, u_1 \ldots)$, $x = (\ldots, x_{-1}, x_0, x_1 \ldots)$. We define the \textit{energy} of a disturbance $w$ to be $$\|w\|_2^2 = \sum_{t = -\infty}^{\infty} \|w_t\|_2^2.$$ As in the finite-horizon setting, we formulate control as an optimization problem, where the goal is select the control actions so as to minimize the quadratic cost 
\begin{equation} \label{ih-lqr-cost}
\sum_{t = -\infty}^{\infty} \left( x_t^*Q x_t + u_t^*R u_t \right),
\end{equation}
where $Q \succeq 0, R \succ 0$; as in the finite-horizon setting, we assume without loss of generality that the system is parameterized so that $R = I$. We assume $\{A, B_u, B_w, Q\}$ are known in advance, so the only uncertainty in the evolution of the system comes from the external disturbance $w$.

We distinguish between several different kinds of information patterns that may be available to a controller. We say a controller is \textit{causal} if in each timestep it is able to observe all previous disturbances up to and including the current timestep, e.g. $u_t = \pi(w_0, \ldots, w_t)$ for some function $\pi$. Similarly, a controller is \textit{strictly causal} if in each timestep it is able to observe all previous disturbances up to but not including the current timestep, e.g. $u_t = \pi(w_0, \ldots, w_{t-1})$. We often use the term \textit{online} to describe causal or strictly causal controllers. A controller is noncausal if it is not causal; in particular, the \textit{clairvoyant offline optimal controller} (sometimes called the noncausal controller) selects the control actions in each timestep with access to the full disturbance sequence $w$ so as to minimize the cost (\ref{fh-lqr-cost}), in the finite-horizon setting, or (\ref{ih-lqr-cost}), in the infinite-horizon setting.

As is standard in the input-output approach to control, we encode controllers as linear \textit{transfer operators} mapping the disturbances to the quadratic cost we wish to minimize. Define $s_t = Q_t^{1/2}x_t$. With this notation, the quadratic costs (\ref{fh-lqr-cost}) and (\ref{ih-lqr-cost}) can be written in a very simple form: $$ \|s\|_2^2 + \|u\|_2^2.$$ The dynamics (\ref{evolution-eq}) are captured by the relation $$s = Fu + Gw, $$ where $F$ and $G$ are strictly causal operators encoding $\{A_t, B_{u, t}, B_{w, t}, Q_t^{1/2}\}_{t=0}^{T-1}$ in the finite-horizon setting and $\{A, B_u, B_w, Q^{1/2}\}$ in the infinite-horizon setting. We refer the reader to \cite{hassibi1999indefinite} for more background on transfer operators and the input-output approach to control.

\subsection{Competitive Control}
The central focus of this paper is designing a controller with optimal competitive ratio:

\begin{problem}[Competitive control] \label{competitive-control-problem}
Find an online controller which minimizes the competitive ratio $$\sup_{w} \frac{ALG(w)}{OPT(w)},$$ where $ALG(w)$ is the cost incurred by the online controller in response to the disturbance $w$ and $OPT(w)$ is the cost incurred by the clairvoyant offline optimal controller. 
\end{problem}

This problem can be studied in both finite-horizon setting and infinite-horizon setting; in the infinite-horizon setting we assume $w$ has bounded energy.
The offline optimal controller has a well-known description at the level of transfer operators (Theorem 11.2.1 in \cite{hassibi1999indefinite}):
\begin{equation}  \label{offline-operator}
u^* = - (I + F^*F)^{-1}F^*Gw.
\end{equation}
Similarly, the offline optimal cost is 
\begin{equation} \label{offline-operator-cost}
OPT(w) = w^*G^*(I + FF^*)^{-1}Gw.
\end{equation}
We note that a state-space description of the offline optimal controller was recently obtained in \cite{goel2020power}.

We call the controller with the smallest possible competitive ratio the \textit{competitive controller}. Instead of minimizing the competitive ratio directly, we instead solve the following relaxation:

\begin{problem}[Suboptimal competitive control] \label{cr-suboptimal-control-problem}
Given $\gamma > 0$, find an online controller such that $$\sup_{w} \frac{ALG(w)}{OPT(w)} < \gamma^2$$ for all disturbances $w$, or determine whether no such controller exists. 
\end{problem}

We call such a controller the \textit{competitive controller at level $\gamma$}. It is clear that if we can solve this suboptimal problem then we can easily recover the competitive controller via bisection on $\gamma$.

\subsection{Robust control} \label{robust-estimation-control-sec}
Our results rely heavily on techniques from robust control. In particular, we show that the problem of obtaining the competitive controller can be reduced to an $H_{\infty}$ control problem:

\begin{problem}[$H_{\infty}$-optimal control] \label{hinf-optimal-control-problem}
Find an online controller that minimizes $$\sup_{w} \frac{ALG(w)}{\|w\|_2^2}, $$ where $ALG(w)$ is the cost incurred by the online controller in response to the disturbance $w$. 
\end{problem}
This problem can be studied in both finite-horizon setting and infinite-horizon setting; in the infinite-horizon setting we assume $w$ has bounded energy.
The $H_{\infty}$-optimal control problem has the natural interpretation of minimizing the worst-case gain from the energy in the disturbance $w$ to the cost incurred by the controller. In general, it is not known how to derive a closed-form for the $H_{\infty}$-optimal controller, so instead is it common to consider a relaxation: 

\begin{problem} [Suboptimal $H_{\infty}$ control at level $\gamma$] \label{hinf-suboptimal-control-problem}
Given  $\gamma > 0$, find an online controller such that  $$ALG(w) < \gamma^2 \| w\|_2^2 $$ for all disturbances $w$, or determine whether no such controller exists. 
\end{problem}
We call such a controller the \textit{$H_{\infty}$ controller at level $\gamma$}. It is clear that if we can solve this suboptimal problem then we can easily recover the $H_{\infty}$-optimal controller via bisection on $\gamma$.
The finite-horizon $H_{\infty}$ controller at level $\gamma$ has a well-known state-space description:
\begin{theorem}[Theorems 9.5.1 and 9.5.2 in \cite{hassibi1999indefinite}] \label{hinf-suboptimal-controller-thm-fh}
Given $\gamma > 0$, a causal finite-horizon $H_{\infty}$ controller at level $\gamma$ exists if and only if
\begin{equation*}
B_{w, t}^*\left[P_{t+1} - P_{t+1}B_{u, t} H_t^{-1}B_{u, t}^*P_{t+1}\right]B_{w, t} \prec \gamma^2I 
\end{equation*}
for all $t = 0, \ldots T-1$, where we define $$H_t = (I + B_{u, t}^*P_{t+1}B_{u, t})$$ and $P_t$ is the solution of the backwards-time Riccati recurrence
\begin{equation*}
P_t = Q_t + A_t^*P_{t+1}A_t -  A_t^*P_{t+1}\tilde{B}_t \tilde{H}_t^{-1} \tilde{B}_t^* P_{t+1}A_t,
\end{equation*}
where we initialize $P_T = 0$, and we define $$ \tilde{B}_t = \begin{bmatrix} B_{u, t} & B_{w, t} \end{bmatrix}, \hspace{3mm} \tilde{R} = \begin{bmatrix} I & 0 \\ 0 & -\gamma^2 I \end{bmatrix},$$ $$\tilde{H}_t =  \tilde{R}  + \tilde{B}_t^* P_{t+1} \tilde{B}_t. $$
In this case, one possible causal finite-horizon $H_{\infty}$ controller at level $\gamma$ is given by
\begin{equation*}
u_t =  -H_t^{-1} B_{u, t}^*P_{t+1}(A_tx_t +  B_{w, t}  w_t).
\end{equation*}
A strictly causal finite-horizon controller at level $\gamma$ exists if and only if $$B_{u, t}^*P_{t+1}B_{u, t} \prec \gamma^2 I$$ for $t = 0 \ldots T - 1$. In this case, one possible strictly causal finite-horizon controller at level $\gamma$ is given by \begin{equation*}
u_t =  -H_t^{-1} B_{u, t}^*P_{t+1}A_tx_t.
\end{equation*}
\end{theorem}

The infinite-horizon $H_{\infty}$ controller at level $\gamma$ also has a well-known state-space description:
\begin{theorem}[Theorem 13.3.3 in \cite{hassibi1999indefinite}] \label{hinf-ih-controller-thm}
Suppose $(A, B_u)$ is stabilizable and $(A, Q^{1/2})$ is observable on the unit circle. A causal controller at level $\gamma$ exists if and only if there exists a solution to the Ricatti equation
\begin{equation*}
P = Q + A^*PA -A^* P \tilde{B}\tilde{H}^{-1}\tilde{B}^*PA
\end{equation*}
with $$ \tilde{B} = \begin{bmatrix} B_{u} & B_{w} \end{bmatrix}, \hspace{3mm} \tilde{R} = \begin{bmatrix} I & 0 \\ 0 & -\gamma^2 I \end{bmatrix},$$ $$\tilde{H} =  \tilde{R}  + \tilde{B}^* P \tilde{B}, $$
such that 
\begin{enumerate}
    \item $ A - \tilde{B}\tilde{H}^{-1}\tilde{B}^*PA$ is stable;
    \item $\tilde{R}$ and $\tilde{H}$ have the same inertia;
    \item $P \succeq 0$.
\end{enumerate}
In this case, the infinite-horizon $H_{\infty}$ controller at level $\gamma$ has the form
\begin{equation*}
u_t =  -H^{-1} B_{u}^*P(Ax_t +  B_{w}  w_t),
\end{equation*}
where $H = I + B_u^*PB_u$. A strictly causal $H_{\infty}$ controller at level $\gamma$ exists if and only if conditions 1 and 3 hold, and additionally $$B_u^*PB_u \prec \gamma^2 I$$ and $$I + B_w^*P(I - \gamma^2 B_u B_u^* P)^{-1}B_u \succ 0.$$ In this case, one possible strictly causal $H_{\infty}$ controller at level $\gamma$ is given by \begin{equation*}
u_t =  -H^{-1} B_{u}^*PAx_t.
\end{equation*}
\end{theorem}

\section{The Competitive Controller}
In this section we present our main results: a computationally efficient state-space description of the competitive controller, i.e. the online controller with the smallest possible competitive ratio, in both the finite-horizon setting and the infinite-horizon setting. In both settings, the key technique we employ is a reduction from the competitive control problem (Problem \ref{competitive-control-problem}) to an $H_{\infty}$ control problem (Problem \ref{hinf-optimal-control-problem}). To perform this reduction, we construct a synthetic dynamical system whose dimension is twice that of the original system. We also construct a new synthetic disturbance $w'$ which can be computed online as the disturbance $w = w_0, w_1, \ldots $ is observed. The $H_{\infty}$ controller in our synthetic system, when fed the synthetic disturbance $w'$, selects the control actions which minimize competitive ratio in the original system. As is standard in $H_{\infty}$ control, we first synthesize the suboptimal $H_{\infty}$ controller
at level $\gamma$; by the nature of our construction, this controller is guaranteed to have competitive ratio at most $\gamma^2$ in the original system. We can then obtain the $H_{\infty}$-optimal controller in the synthetic system (and hence the competitive controller in the original system) by minimizing $\gamma$ subject to the constraints outlined in Theorems \ref{hinf-suboptimal-controller-thm-fh} and \ref{hinf-ih-controller-thm}.

Recall that $$OPT(w) = w^*G^*(I + FF^*)^{-1}Gw.$$
It follows that Problem \ref{cr-suboptimal-control-problem} can be expressed as finding an online  controller such that 
\begin{equation} \label{cr-suboptimal-cond}
ALG(w) < \gamma^2w^*G^*(I + FF^*)^{-1}Gw
\end{equation}
for all disturbances $w$, or determining whether no such controller exists. Let $\Delta$ be the unique casual operator such that $\Delta \Delta^* = I + FF^*$. Then condition (\ref{cr-suboptimal-cond}) can be rewritten as an $H_{\infty}$ condition:
$$ALG(w) < \gamma^2 \|w'\|_2^2, $$ 
where we define $w' = \Delta^{-1}Gw$. The dynamical system $s = Fu + Gw$, which is driven by the disturbance $w$, can be transformed into a system driven by $w'$:  
\begin{eqnarray}
s &=& Fu + Gw \nonumber \\
&=& Fu + (\Delta \Delta^{-1})Gw \nonumber \\
&=&  Fu + \Delta  (\Delta^{-1}Gw) \nonumber \\
&=&  Fu + \Delta w'. \label{transformed-dynamics}
\end{eqnarray}
We have shown that the problem of finding a competitive controller at level $\gamma$ in the system $s = Fu + Gw$ is equivalent to finding an $H_{\infty}$ controller at level $\gamma$ in the system $s = Fu + \Delta w'$; the key is to obtain the factorization $\Delta \Delta^* = I + FF^*$. In the finite-horizon setting, we obtain this factorization using state-space models and the whitening property of the Kalman filter; in the infinite-horizon setting we first pass to the frequency domain and employ algebraic techniques to factor $\Delta(z) \Delta^*(z^{-*}) = I + F(z)F^*(z^{-*})$, and then reconstruct the controller in time domain from its frequency domain model.

\subsection{Finite-horizon competitive control}
We first consider finite-horizon control in linear time-varying systems as described in Section \ref{preliminaries-sec}. We prove:

\begin{theorem}[Finite-horizon competitive control] \label{competitive-controller-thm-fh}
A causal finite-horizon controller with competitive ratio bounded above by $\gamma^2$ exists if and only if 
\begin{equation} \label{competitive-gamma-condition}
 \hat{B}_{w, t}^*\left[\hat{P}_{t+1} - \hat{P}_{t+1}\hat{B}_{u, t} \hat{H}_t^{-1}\hat{B}_{u, t}^*{P}_{t+1}\right]\hat{B}_{w, t} \prec \gamma^2I 
\end{equation}
for $t = 0, \ldots, T-1$, where we define $$\hat{A}_t = \begin{bmatrix} A_t & K_t \Sigma_t^{1/2} \\ 0 & 0 \end{bmatrix}, \hspace{3mm} \hat{B}_{u, t} = \begin{bmatrix} B_{u, t} \\ 0 \end{bmatrix}, \hspace{3mm} \hat{B}_{w, t} = \begin{bmatrix} 0 \\ I \end{bmatrix}, $$ $$ \hat{Q}_t = \begin{bmatrix} Q_t & Q_t^{1/2} \Sigma_t^{1/2} \\ \Sigma_t^{1/2}Q_t^{1/2} & \Sigma_t \end{bmatrix}, \hspace{3mm} \hat{H}_t = I + \hat{B}_{u, t}^*\hat{P}_{t+1} \hat{B}_{u, t}. $$ 
we define $\hat{P}_t$ to be the solution of the backwards-time Riccati recursion 
\begin{equation} \label{competitive-fh-ricatti-recur}
\hat{P}_t = \hat{Q}_t + \hat{A}_t^*\hat{P}_{t+1}\hat{A}_t -  \hat{A}_t^*\hat{P}_{t+1}\tilde{B}_t \tilde{H}^{-1} \tilde{B}_t^* \hat{P}_{t+1}\hat{A}_t 
\end{equation}
where we initialize $\hat{P}_T = 0$ and define $$ \tilde{B}_t = \begin{bmatrix} \hat{B}_{u, t} & \hat{B}_{w, t} \end{bmatrix},$$ $$\tilde{H}_t = \begin{bmatrix} I & 0 \\ 0 & -\gamma^2 I \end{bmatrix}  + \tilde{B}_t^* \hat{P}_{t+1} \tilde{B}_t,$$ and $K_t, \Sigma_t$ are defined in (\ref{K-and-Sigma}).
In this case, a causal controller with competitive ratio bounded above by $\gamma^2$ is given by $$ u_t =  -\hat{H}_t^{-1} \hat{B}_{u, t}^*\hat{P}_{t+1}\left(\hat{A}_t \xi_t +  \hat{B}_{w, t}  w_{t+1}' \right),$$
where the dynamics of $\xi$ are
\begin{equation} \label{final-ss}
\xi_{t+1} = \hat{A}_t\xi_t + \hat{B}_{u, t} u_t + \hat{B}_{w, t} w_{t+1}'
\end{equation}
and we initialize $\xi_0 = 0$. The synthetic disturbance $w'$ can be computed using the recursion
$$\nu_{t+1} = (A_t - K_t Q_t^{1/2})\nu_t + B_{w, t}w_t, \hspace{3mm} w_t' = \Sigma_t^{-1/2} Q_t^{1/2} \nu_t,$$ where we initialize $\nu_0 = 0$.  A strictly causal finite-horizon controller with competitive ratio bounded above by $\gamma^2$ exists if and only if 
\begin{equation*}
 \hat{B}_{w, t}^*\hat{P}_{t+1} \hat{B}_{w, t} \prec \gamma^2I 
\end{equation*}
for $t = 0, \ldots, T-1$. In this case, a strictly causal controller with competitive ratio bounded above by $\gamma^2$ is given by $$ u_t =  -\hat{H}_t^{-1} \hat{B}_{u, t}^*\hat{P}_{t+1}\hat{A}_t \xi_t.$$
\end{theorem}

We make a few observations. First, comparing with Theorem \ref{hinf-suboptimal-controller-thm-fh}, we see that the competitive controller at level $\gamma$ has a similar structure to the $H_{\infty}$ controller at level $\gamma$; indeed, the competitive controller is just the $H_{\infty}$ controller in the system (\ref{final-ss}). Second, we emphasize that the synthetic disturbance $w'$ appearing in (\ref{final-ss}) is a \textit{strictly causal} function of $w$; in particular, $w_{t+1}'$ is a linear combination of the disturbances $w_0, \ldots, w_t$. This is crucial, since it means that we can construct $w_{t+1}'$ online, using only the observations available up to time $t$. Third, we note that the dimension of the control input $u$ in the synthetic system (\ref{final-ss}) is the same as the dimension of the control input in the original system (\ref{evolution-eq}); this allows us to use $u$ to steer the original system. Lastly, since the competitive controller is simply the standard $H_{\infty}$ controller in a system of dimension $2n$, it is clear that the computational resources required to implement the competitive controller are identical to those required to implement the $H_{\infty}$ controller, up to a constant factor.

The proof of Theorem \ref{competitive-controller-thm-fh} is presented in the appendix.

\subsection{Infinite-horizon competitive control}
We next consider  infinite-horizon control in linear time-invariant systems as described in Section \ref{preliminaries-sec}. We prove:

\begin{theorem}[Infinite-horizon competitive control] \label{competitive-controller-thm-ih}
Suppose $(A, B_u)$ is stabilizable and $(A, Q^{1/2})$ is detectable. A causal infinite-horizon controller with competitive ratio bounded above by $\gamma^2$ exists if and only if there exists a solution to the Ricatti equation
\begin{equation} \label{competitive-ih-ricatti-eq}
\hat{P} = \hat{Q} + \hat{A}^*\hat{P}\hat{A} -\hat{A}^* \hat{P} \tilde{B}\tilde{H}^{-1}\tilde{B}^*\hat{P}\hat{A}
\end{equation}
with $$\hat{A} = \begin{bmatrix} A & K \Sigma^{1/2} \\ 0 & 0 \end{bmatrix}, \hspace{3mm} \hat{B}_{u} = \begin{bmatrix} B_{u} \\ 0 \end{bmatrix}, \hspace{3mm} \hat{B}_{w} = \begin{bmatrix} 0 \\ I \end{bmatrix}, $$ $$ \hat{Q} = \begin{bmatrix} Q & Q^{1/2} \Sigma^{1/2} \\ \Sigma^{1/2}Q^{1/2} & \Sigma \end{bmatrix},$$ $$\tilde{B} = \begin{bmatrix} \hat{B}_{u} & \hat{B}_{w} \end{bmatrix}, \hspace{3mm} \tilde{R} = \begin{bmatrix} I & 0 \\ 0 & -\gamma^2 I \end{bmatrix},$$ $$\tilde{H} =  \tilde{R}  + \tilde{B}^* P \tilde{B}, $$
and $K, \Sigma$ defined in (\ref{sigma-k-ih-definition}), such that 
\begin{enumerate}
    \item $ \hat{A} - \tilde{B}\tilde{H}^{-1}\tilde{B}^*\hat{P}\hat{A}$ is stable;
    \item $\tilde{R}$ and $\tilde{H}$ have the same inertia;
    \item $\hat{P} \succeq 0$.
\end{enumerate}
In this case, a causal infinite-horizon $H_{\infty}$ controller at level $\gamma$ is given by
\begin{equation*}
u_t =  -\hat{H}^{-1} \hat{B}_{u}^*\hat{P}(\hat{A}\xi_t +  \hat{B}_w  w_{t+1}'),
\end{equation*}
where $\hat{H} = I + \hat{B}_u^*\hat{P}\hat{B}_u$ and the dynamics of $\xi$ are
\begin{equation*}
\xi_{t+1} = \hat{A} \xi_t + \hat{B}_{u} u_t + \hat{B}_{w} w_{t+1}' 
\end{equation*}
and we initialize $\xi_0 = 0$. The synthetic disturbance $w'$ can be computed using the recursion
$$\nu_{t+1} = (A - K Q^{1/2})\nu_t + B_{w}w_t, \hspace{3mm} w_t' = \Sigma^{-1/2} Q^{1/2} \nu_t,$$ where we initialize $\nu_0 = 0$. A strictly causal infinite-horizon controller with competitive ratio bounded above by $\gamma^2$ exists if and only if conditions 1 and 3 hold, and additionally $$\hat{B}_u^*\hat{P}\hat{B}_u \prec \gamma^2 I$$ and $$I + \hat{B}_w^*\hat{P}(I - \gamma^2 \hat{B}_u \hat{B}_u^* \hat{P})^{-1}\hat{B}_u \succ 0.$$
In this case, a strictly causal controller with competitive ratio bounded above by $\gamma^2$ is given by \begin{equation*}
u_t =  -\hat{H}^{-1} \hat{B}_{u}^*\hat{P}\hat{A}\xi_t.
\end{equation*}
\end{theorem}

We note that the infinite-horizon controller described in Theorem \ref{competitive-controller-thm-ih} is identical to the finite-horizon controller described in Theorem \ref{competitive-controller-thm-fh}, except that the Ricatti recursion (\ref{competitive-fh-ricatti-recur}) is replaced by a Ricatti equation (\ref{competitive-ih-ricatti-eq}), and all the matrices appearing in the controller are time-invariant; this is consistent with our intuition that the infinite-horizon controller is the finite-horizon controller in steady-state, in the asymptotic limit as the time-horizon $T$ tends to infinity. It is clear that an infinite-horizon competitive controller with competitive ratio bounded by $\gamma$ is stabilizing (whenever such a controller exists), because its cost is always at most a factor of $\gamma^2$ more than the offline optimal cost, and the offline controller is stabilizing.

The proof of Theorem \ref{competitive-controller-thm-ih} is presented in the appendix.

\section{Numerical Experiments}
We benchmark the causal infinite-horizon competitive controller against the $H_2$-optimal, $H_{\infty}$-optimal, and offline optimal controllers in both a linear system and a nonlinear system. 

\subsection{Boeing 747 Flight Control}

We consider the longitudinal flight control system of a Boeing 747 with linearized dynamics. Assuming level flight at 40,000ft at a speed of 774ft/sec and a discretization interval of 1 second, the dynamics are given by 
\[x_{t+1} = Ax_t + B u_t + w_t, \]
where

$$ A = \begin{bmatrix}0.99 & 0.03 & -0.02 & -0.32 \\ 0.01 & 0.47 & 4.7 & 0.0 \\ 0.02 & -0.06& 0.40& 0.0 \\ 0.01 & -0.04 & 0.72 & 0.99\end{bmatrix}, $$
$$ B = \begin{bmatrix} 0.01 & 0.99 \\ -3.44 & 1.66 \\ -0.83 & 0.44 \\ -0.47 & 0.25\end{bmatrix}. $$
The state $x$ consists of kinematic variables such as velocity and orientation and the control inputs are thrust and elevator angle; we refer to \cite{boeing747} for more information. We assume the initial condition $x_0 = 0$ and take $Q, R = I$. 

We synthesize the infinite-horizon competitive controller using Theorem \ref{competitive-controller-thm-ih} and find that the smallest choice of $\gamma$ satisfying the constraints (\ref{competitive-gamma-condition}) is $\gamma = 1.33$, so the competitive ratio of the competitive controller is $\gamma^2 = 1.77$. In other words, \textit{the cost incurred by our competitive controller is guaranteed to always be within 77\% of the cost incurred by the clairvoyant offline optimal controller, no matter the input disturbance}. We emphasize that our competitive controller is guaranteed to obtain the smallest possible competitive ratio among all online controllers, therefore no online controller can achieve a competitive ratio less than 1.77.

\begin{figure}[h!]
\centering
\includegraphics[width=\columnwidth]{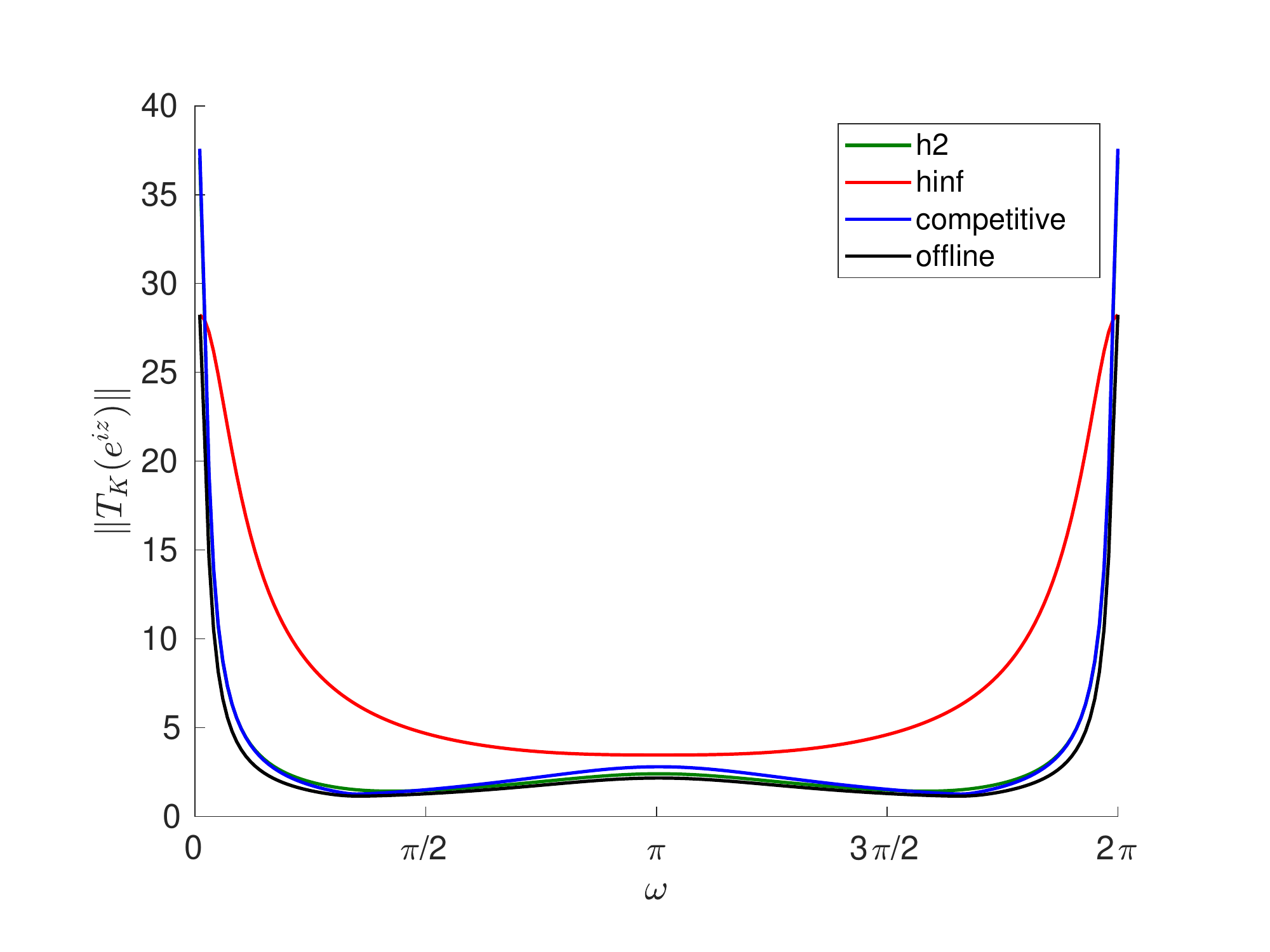}
\caption{Frequency response in the Boeing 747 flight control system. }
\label{boeing747-freq-response-fig}
\end{figure}

\begin{figure}[h!]
\centering
\includegraphics[width=\columnwidth]{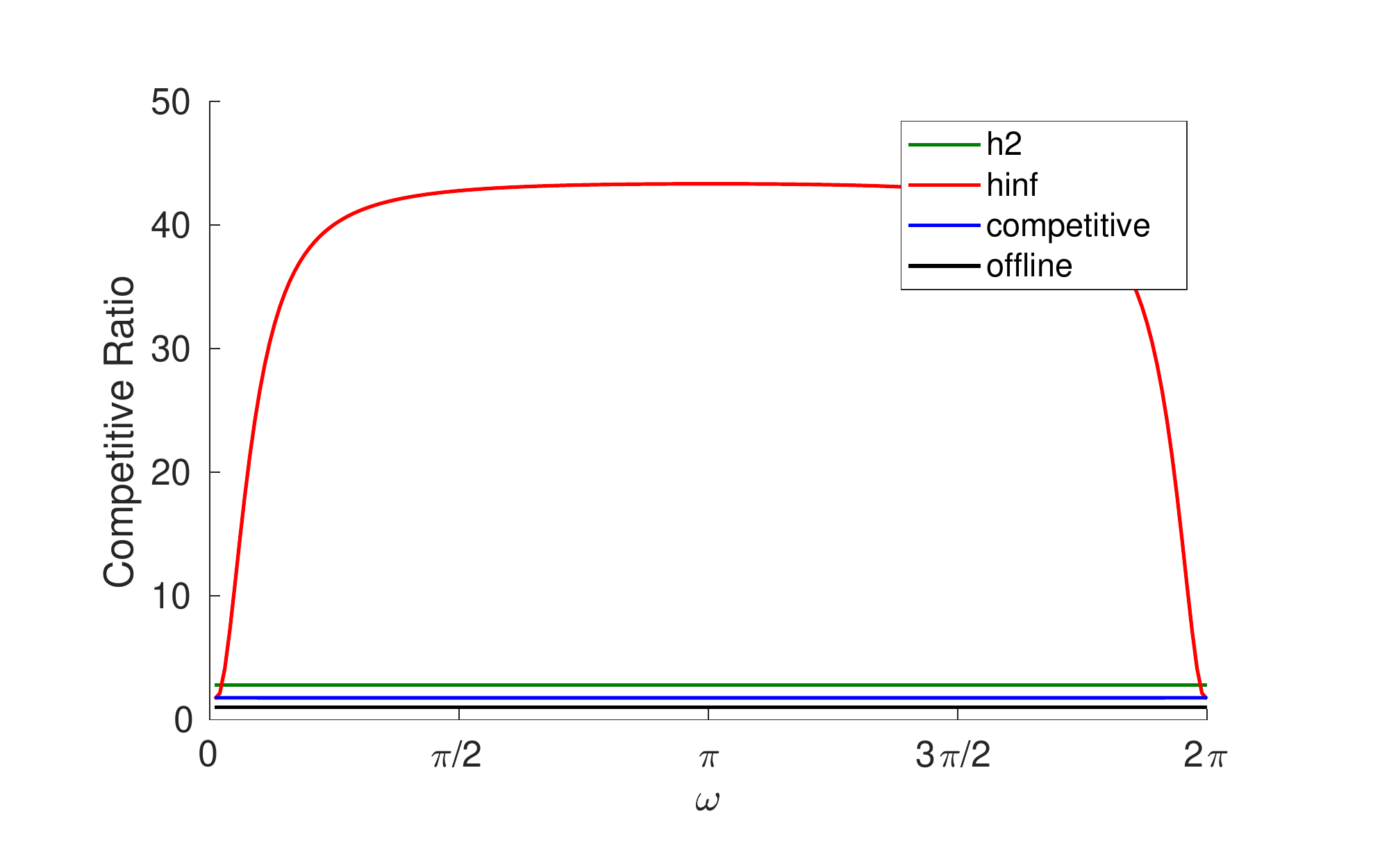}
\caption{Competitive ratio of LQ controllers in a Boeing 747 flight control system across various noise frequencies. }
\label{boeing747-cr-fig}
\end{figure}

In Figure \ref{boeing747-freq-response-fig} we plot the magnitude of $T_K(e^{i\omega})$ at various frequencies $\omega$; this measures how much energy is transferred from the input disturbance to the control cost at the frequency $\omega$. The $H_{\infty}$ controller is designed to be robust to disturbances at all frequencies and hence has the lowest peak. Both the competitive controller and the $H_2$-optimal controller closely track the offline optimal controller.
In Figure \ref{boeing747-cr-fig} we plot the competitive ratio of the various controllers across various frequencies. We see that the competitive ratio of the $H_{\infty}$ controller can be as high as 43.3 at certain frequencies, while the competitive ratio of the $H_2$-optimal controller is 2.8 at every frequency. We note that the competitive ratio of the competitive controller is the smallest at 1.77, as expected. 

\begin{figure}[h!]
\centering
\includegraphics[width=\columnwidth]{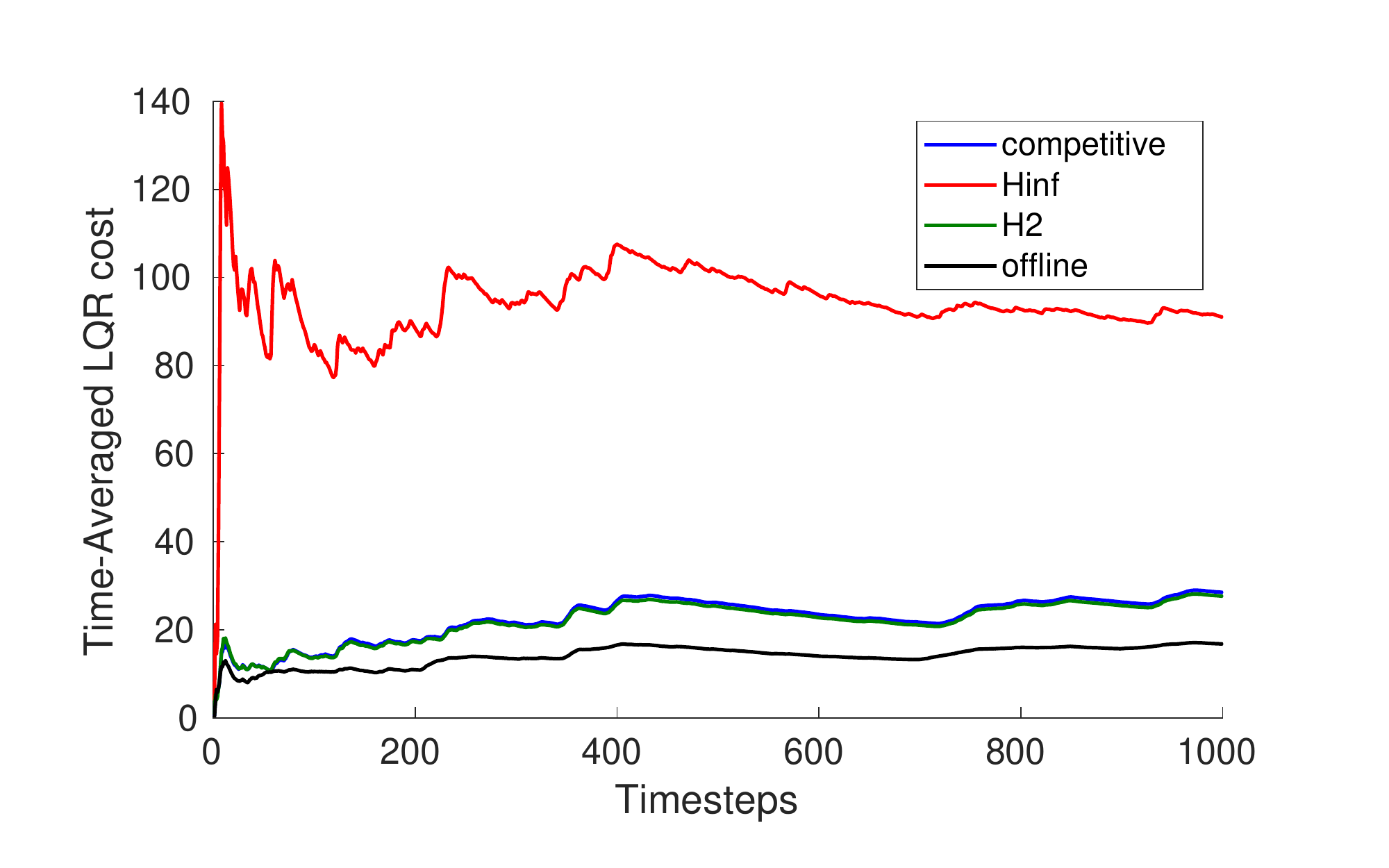}
\caption{Relative performance of LQ controllers in a Boeing 747 flight control system driven by white noise. The competitive controller almost matches the performance of the H2 controller, without being tuned for stochastic noise.}
\label{boeing747-gaussian-fig}
\end{figure}

\begin{figure}[h!]
\centering
\includegraphics[width=\columnwidth]{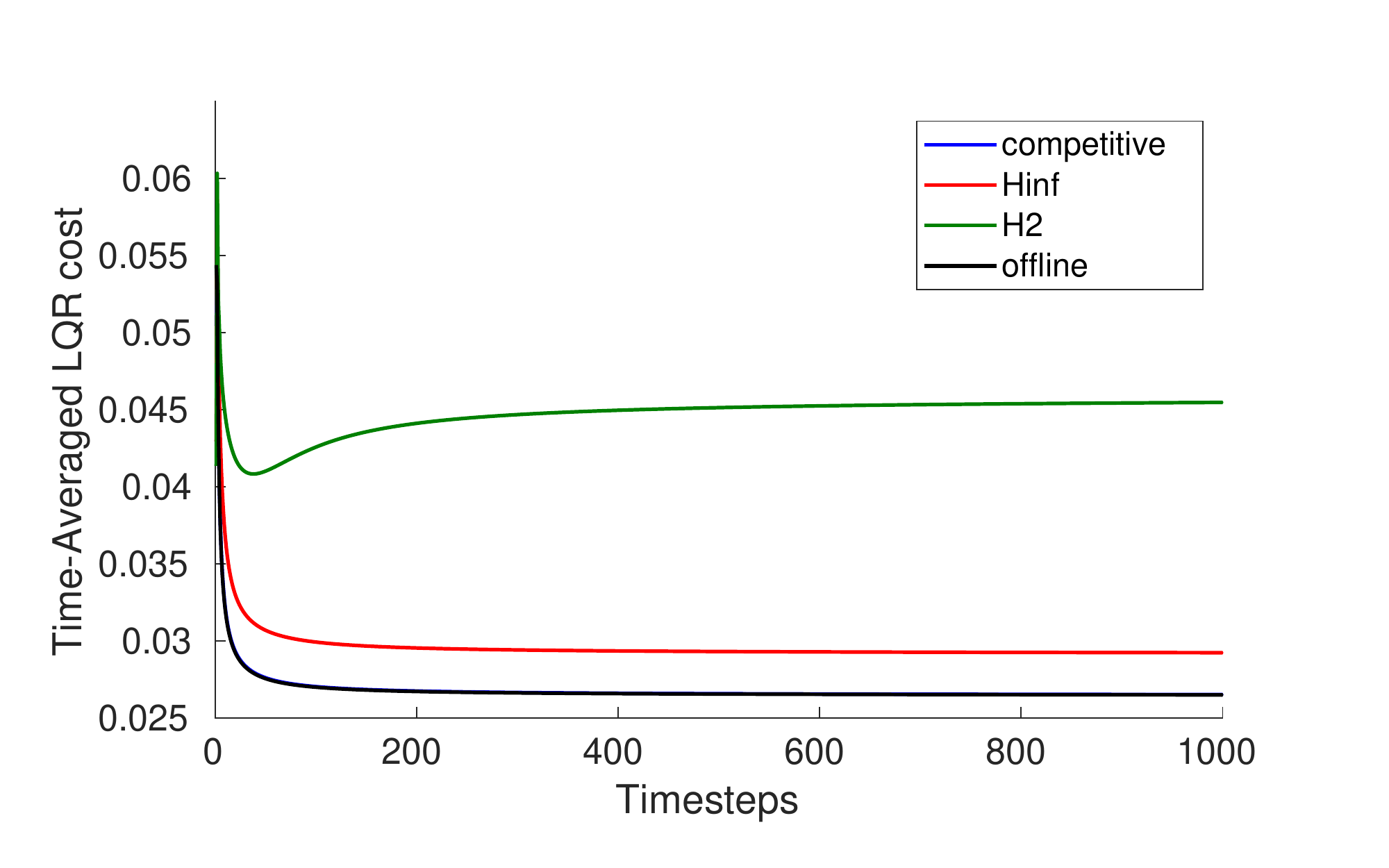}
\caption{Relative performance of LQ controllers in a Boeing 747 flight control system driven by best-case DC noise. The competitive controller exactly matches the performance of the offline optimal controller; the H2 controller incurs substantially more cost.}
\label{boeing747-best-case-fig}
\end{figure}

\begin{figure}[h!]
\centering
\includegraphics[width=\columnwidth]{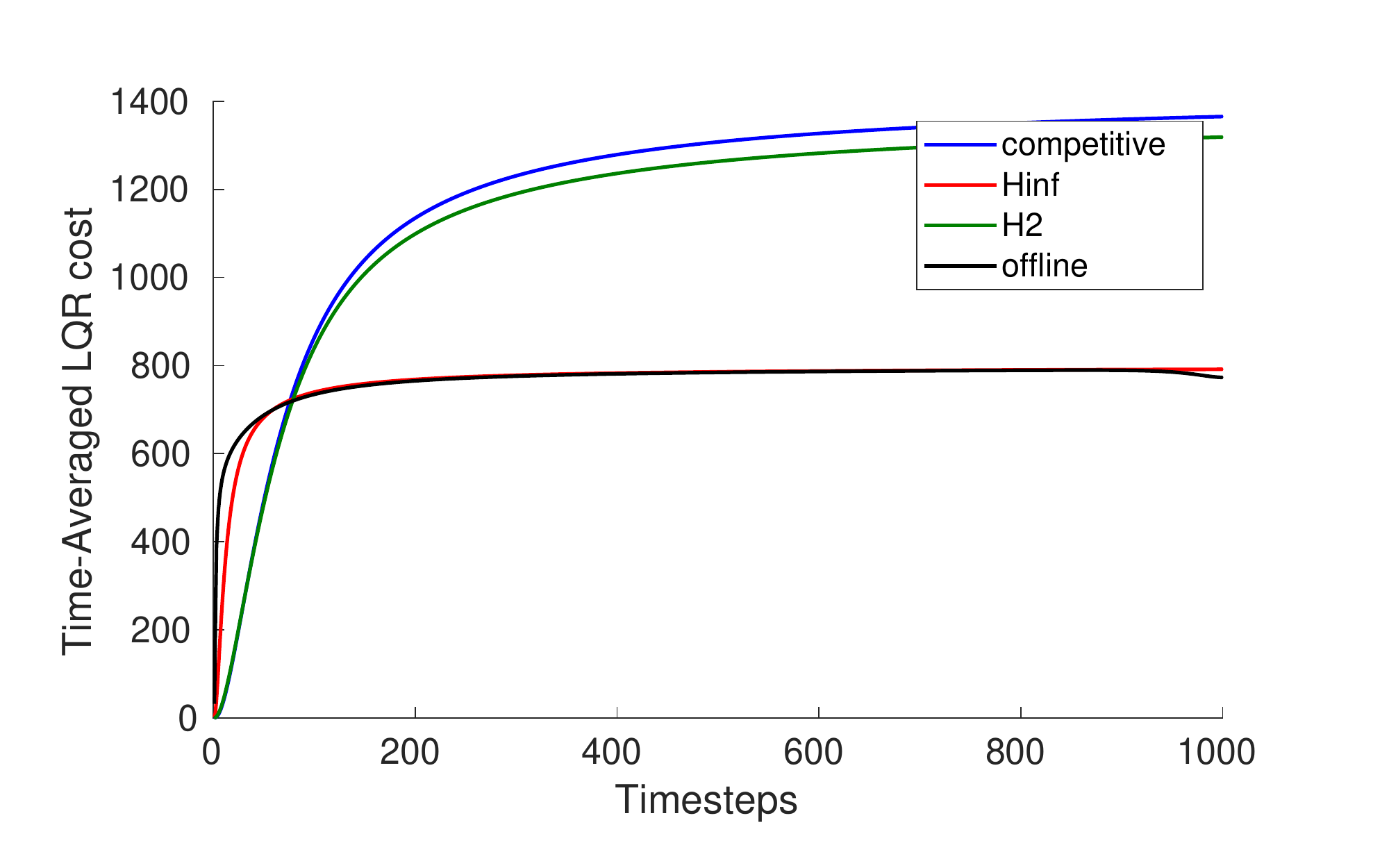}
\caption{Relative performance of LQ controllers in a Boeing 747 flight control system driven by worst-case DC noise. The competitive controller incurs the most cost, but its cost is guaranteed to be at most 77\% more than the cost incurred by the offline optimal controller.}
\label{boeing747-worst-case-fig}
\end{figure}

\begin{figure}[h!]
\centering
\includegraphics[width=\columnwidth]{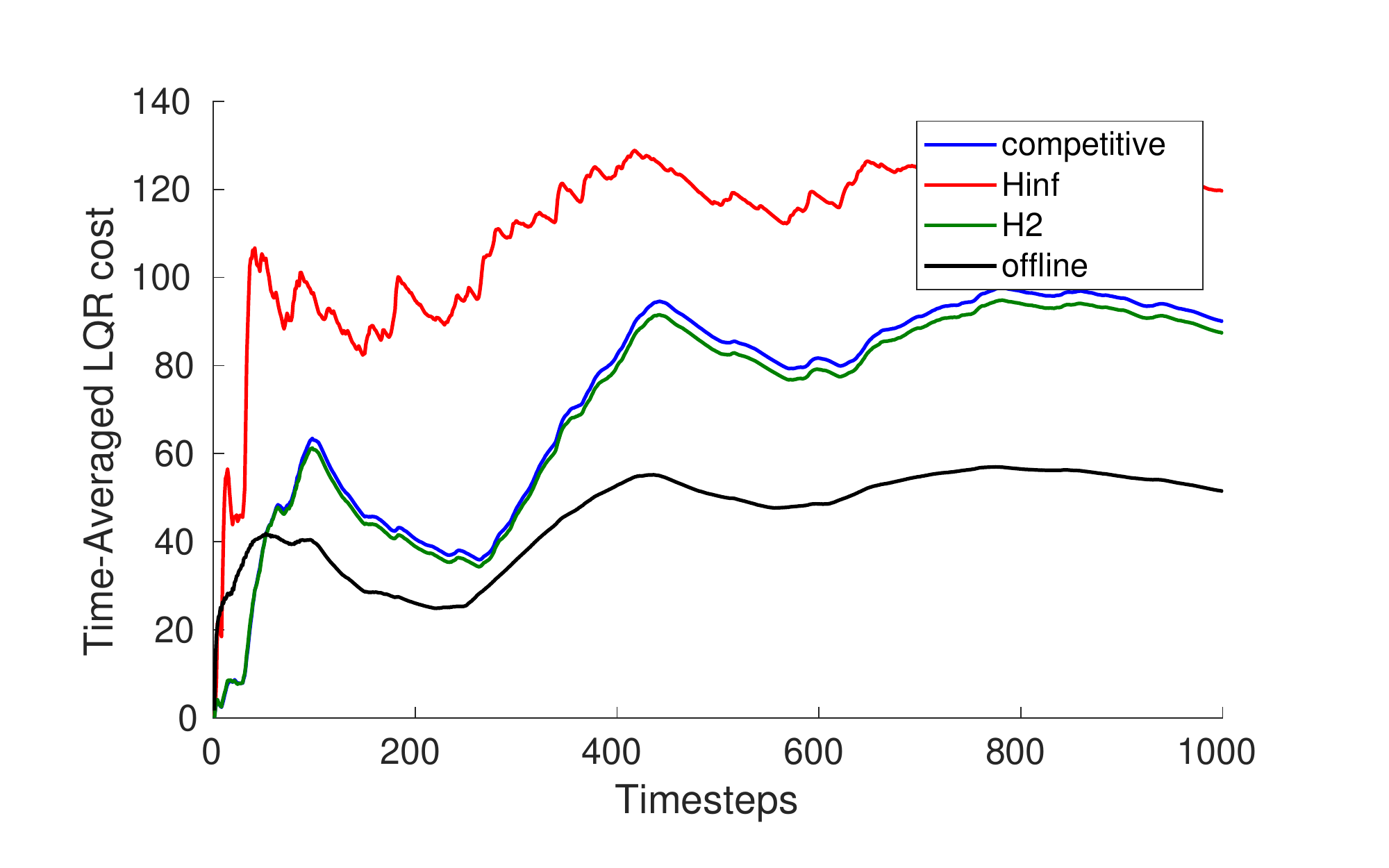}
\caption{Relative performance of LQR controllers in a Boeing 747 flight control system driven by noise which is a mixture of white and worst-case DC components. The competitive controller almost matches the H2 controller and outperforms the Hinf controller.}
\label{boeing747-mixture-fig}
\end{figure}

We next compare the performance of the competitive controller and the $H_2$-optimal, $H_{\infty}$-optimal, and clairvoyant offline optimal controllers across several input disturbances which capture average-case, best-case, and worst-case scenarios for the competitive controller. In Figure \ref{boeing747-gaussian-fig} we plot the controllers' performance when the driving disturbance is white Gaussian noise; unsurprisingly, the $H_2$-optimal controller incurs the lowest cost. The competitive controller is almost able to match the performance of the $H_2$ controller, despite not being designed specifically for stochastic disturbances. We next calculate the best-case and worst-case DC disturbances by computing the eigenvectors corresponding to the smallest and largest eigenvalues of $T_K(e^{i\omega})^*T_K(e^{i\omega})$ at $\omega = 0$, where $T_K$ is the transfer operator associated to the competitive controller. In Figure \ref{boeing747-best-case-fig}, we plot the controller's performance when the noise is taken to be the best-case DC component. The competitive controller exactly matches the performance of the clairvoyant noncausal controller, outperforming the $H_{\infty}$-optimal controller and greatly outperforming the $H_{2}$-optimal controllers. We next plot the controller's performance when the noise is taken to be the worst-case DC component in Figure \ref{boeing747-worst-case-fig}. The competitive controller incurs the highest cost; this is unsurprising, since the noise is chosen specifically to penalize the competitive controller. We note that the ratio of the competitive controller's cumulative cost to that of the offline optimal controller slowly approaches 1.77 as predicted by our competitive ratio bound.
Lastly, in Figure \ref{boeing747-mixture-fig}, we plot the controllers' performance when the noise is a mixture of white and worst-case DC components; we see that the competitive controller almost matches the performance of the $H_2$ controller. Together, these plots highlight the best-of-both-worlds behavior of the competitive controller: in best-case or average-case scenarios it matches or outperforms standard $H_2$ and $H_{\infty}$ controllers, while in the worst-case scenario it is never worse by more than a factor of 1.77.

\subsection{Inverted Pendulum}

We also benchmark our competitive controller in a nonlinear inverted pendulum system. This system has two scalar states, $\theta$ and $\dot{\theta}$, representing angular position and angular velocity, respectively, and a single scalar control input $u$. The state $(\theta, \dot{\theta})$ evolves according to the nonlinear evolution equation 
$$\frac{d}{dt} \begin{bmatrix} \theta \\ \dot{\theta} \end{bmatrix} = \begin{bmatrix} \dot{\theta} \\ \frac{mg\ell}{J}\sin{\theta} + \frac{\ell}{J}u\cos{\theta} + \frac{\ell}{J}w\cos{\theta} \end{bmatrix},$$ 
where $w$ is an external disturbance, and $m, \ell, g, J$ are physical parameters describing the system. Although these dynamics are nonlinear, we can benchmark the regret-optimal controller against the $H_2$-optimal, $H_{\infty}$-optimal, and clairvoyant offline optimal controllers using Model Predictive Control (MPC). In the MPC framework, we iteratively linearize the model dynamics around the current state, compute the optimal control signal in the linearized system, and then update the state in the original nonlinear system using this control signal. In our experiments we take $Q, R = I$ and initialize $\theta$ and $\dot{\theta}$ to zero. We assume that units are scaled so that all physical parameters are 1. We set the discretization parameter $\delta_t = 0.001$ and sample the dynamics at intervals of $\delta_t$.

\begin{figure}[h!]
\centering
\includegraphics[width=\columnwidth]{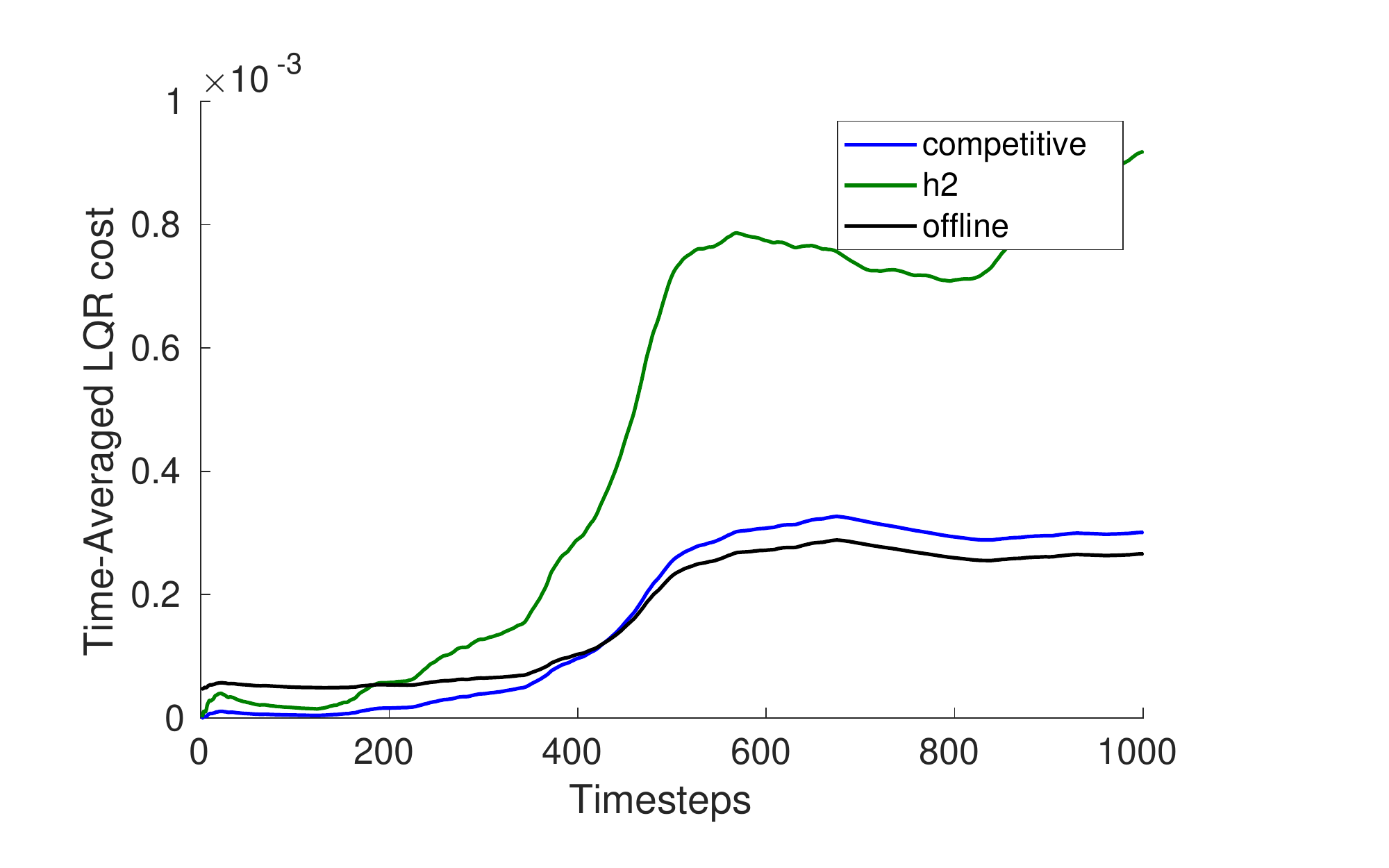}
\caption{Relative performance of LQ controllers in an inverted pendulum system driven by zero-mean Gaussian noise. }
\label{invpendgaussian-fig}
\end{figure}

\begin{figure}[h!]
\centering
\includegraphics[width=\columnwidth]{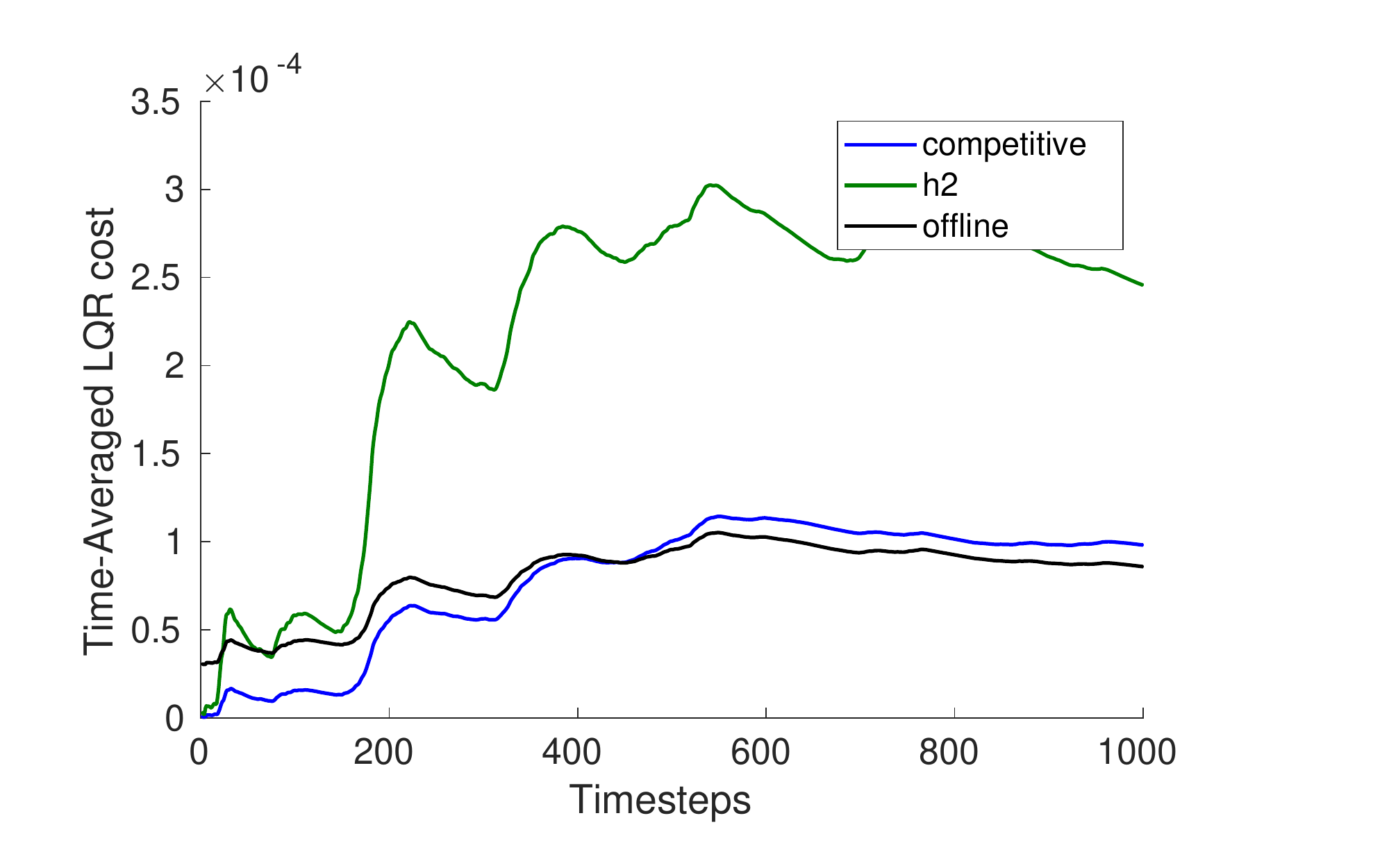}
\caption{Relative performance of LQ controllers in an inverted pendulum system driven by Gaussian noise whose mean varies sinusoidally over time.}
\label{invpendsine_gaussian-fig}
\end{figure}

\begin{figure}[h!]
\centering
\includegraphics[width=\columnwidth]{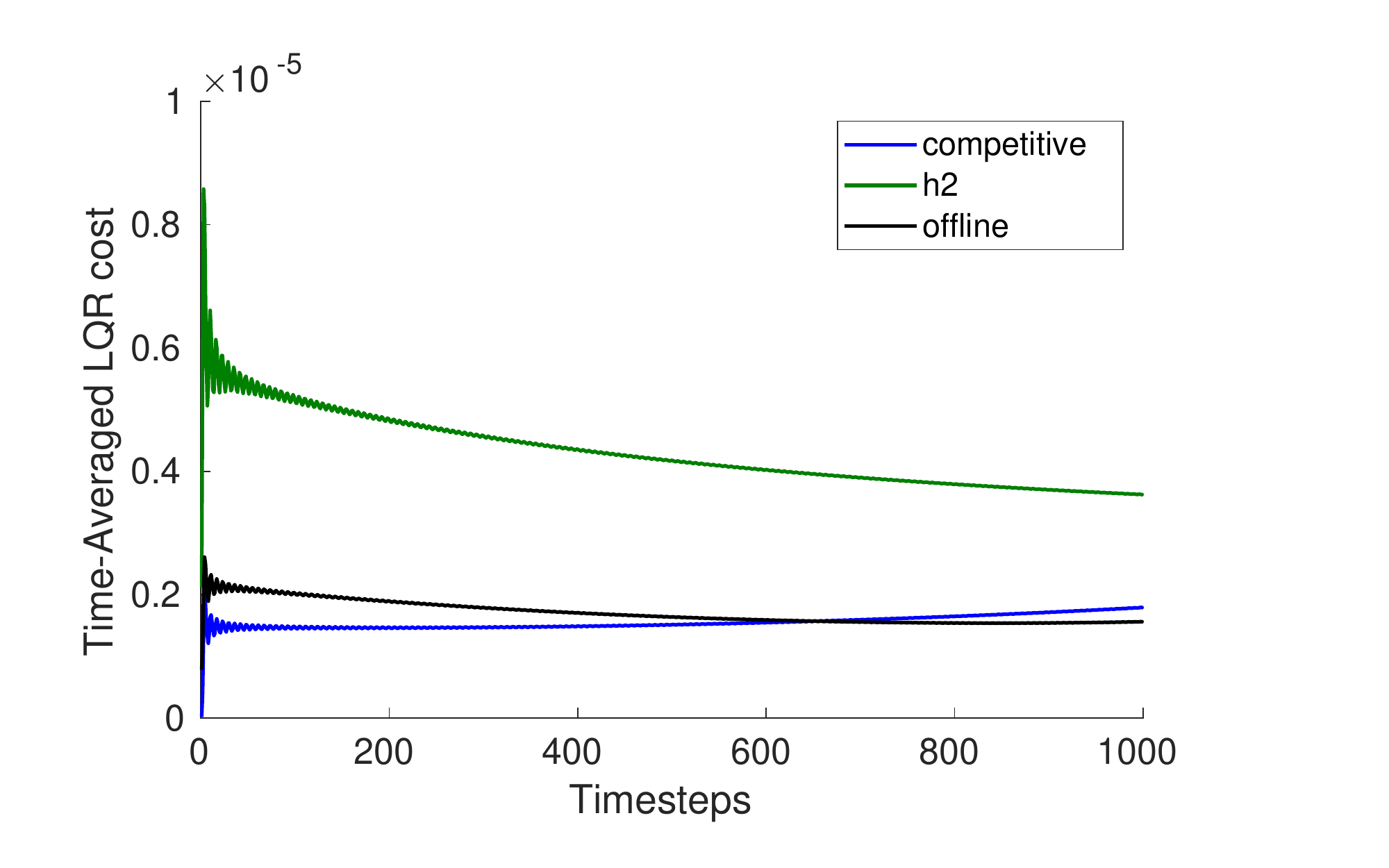}
\caption{Relative performance of LQ controllers in an inverted pendulum system driven by a high-frequency sinusoidal disturbance.}
\label{invpend-high-freq-sine-fig}
\end{figure}

\begin{figure}[h!]
\centering
\includegraphics[width=\columnwidth]{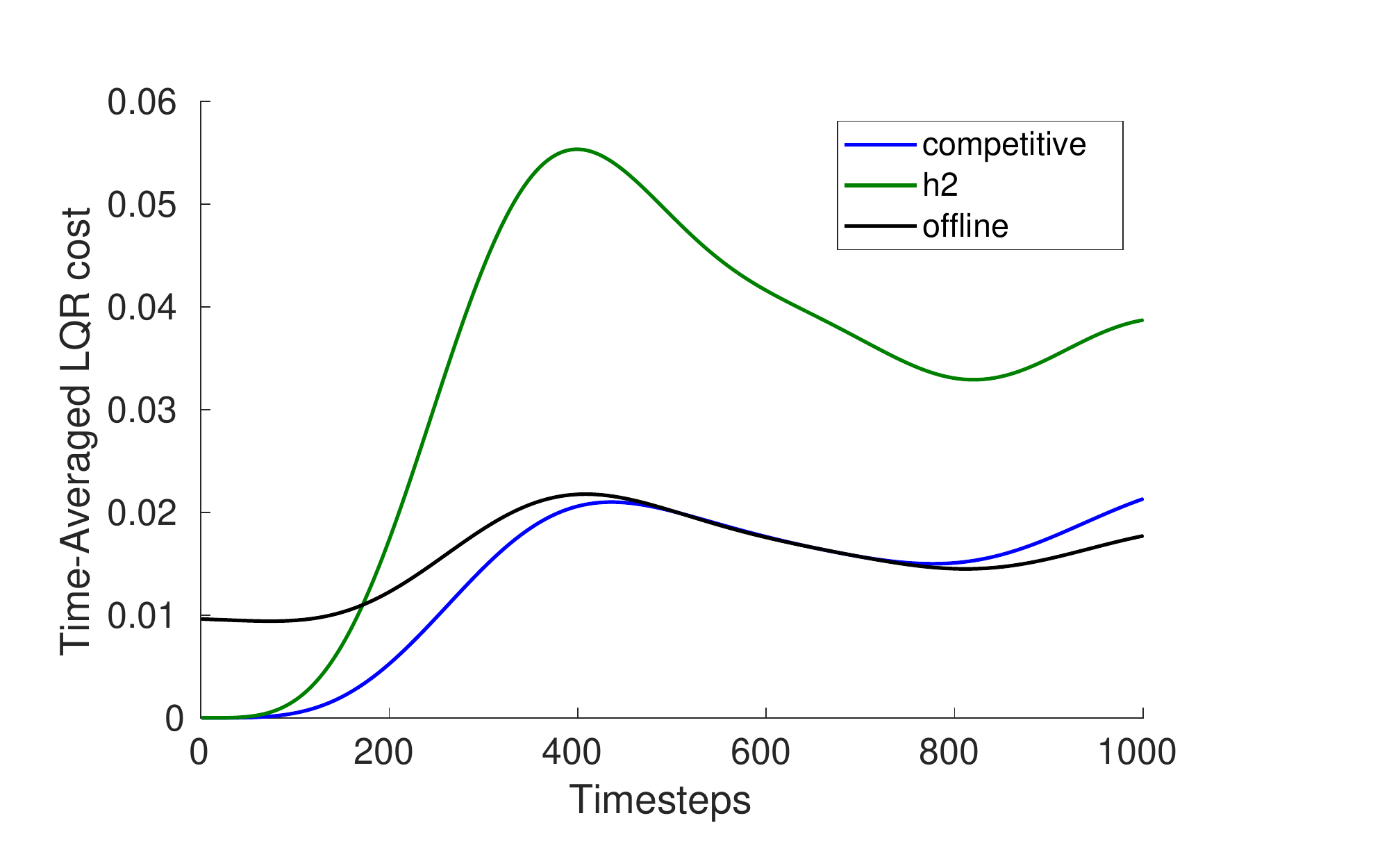}
\caption{Relative performance of LQ controllers in an inverted pendulum system driven by a low-frequency sinusoidal disturbance.}
\label{invpend-low-freq-sine-fig}
\end{figure}

\begin{figure}[h!]
\centering
\includegraphics[width=\columnwidth]{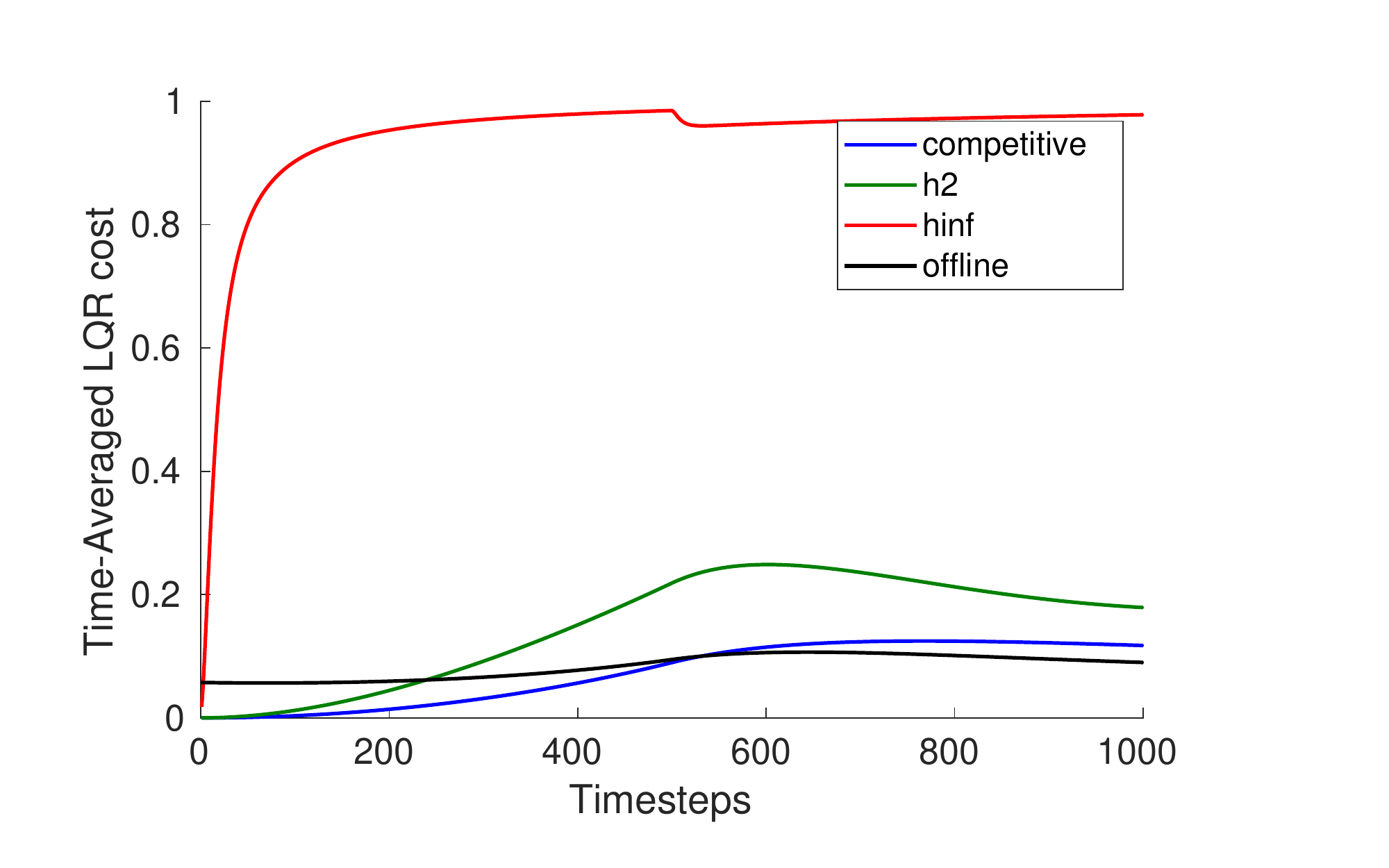}
\caption{Relative performance of LQ controllers in an inverted pendulum system driven by step-function noise.}
\label{invpendstep-fig}
\end{figure}

In Figure \ref{invpendgaussian-fig}, we plot the relative performance of the various controllers when the noise is drawn i.i.d. from a standard Gaussian distribution in each timestep. Surprisingly, the competitive controller significantly outperforms the $H_2$-optimal controller, which is tuned for i.i.d zero-mean noise; this may be because the competitive controller is better able to adapt to nonlinear dynamics.  The cost incurred by the $H_{\infty}$-optimal controller is orders of magnitude larger than that of the other controllers and is not shown. 
In Figur \ref{invpendsine_gaussian-fig}, the noise is drawn from a Gaussian distribution whose variance is fixed but whose mean varies over time; we take $w(k\delta_t) \sim \mathcal{N}(\sin(k\delta), 1)$, for $k = 0, \ldots, 1000$. The $H_2$-optimal controller incurs roughly three times the cost of the competitive controller, while the competitive controller closely tracks the performance of the offline optimal controller. As before, the cost incurred by the $H_{\infty}$-optimal controller is orders of magnitude larger than that of the other controllers and is not shown. In Figures \ref{invpend-high-freq-sine-fig} and \ref{invpend-low-freq-sine-fig} we compare the competitive controller to the $H_2$-optimal and offline optimal controllers with both high frequency and low frequency sinusoidal disturbances, with no Gaussian component, e.g $w(k\delta_t) = \sin(k\delta_t)$ and $w(k\delta_t) = \sin(0.01k\delta_t)$. In both plots the competitive controller easily beats the $H_2$ controller and nearly matches the performance of the offline optimal controller.

Lastly, in Figure \ref{invpendstep-fig} we plot the controllers' performance when the noise is generated by a step-function: for the first 500 timesteps the input disturbance is equal to 1, and for the next 500 timesteps it is equal to $-1$. The sudden transition at $t = 500$ presents a challenge for controllers which adapt online, since the new set of input disturbances is completely different than those which had been observed previously. We see that the competitive-controller closely tracks the offline optimal controller and easily outperforms the $H_2$-optimal and $H_{\infty}$-optimal controllers.

\section{Conclusion}
We introduce a new class of controllers, competitive controllers, which dynamically adapt to the input disturbance so as to track the performance of the clairvoyant offline optimal controller as closely as possible. The key idea is to extend classical $H_{\infty}$ control, which seeks to design online controllers so as to minimize the ratio of their control cost to the energy in the disturbance, to instead minimize competitive ratio. We derive the competitive controller in both finite-horizon, time-varying systems and in infinite-horizon, time-invariant systems. 
In both settings, the key idea is to construct a synthetic system and a synthetic driving disturbance such that the $H_{\infty}$-optimal controller in the synthetic system selects the control outputs which minimize competitive ratio in the original system. The main technical hurdle in our construction is the factorization of certain algebraic expressions involving the transfer operator associated to the offline optimal controller. In the finite-horizon setting, we perform this factorization in time domain using the whitening property of the Kalman filter, whereas in the infinite-horizon setting we perform the factorization in frequency domain and then reconstruct the controller in time domain. 

We benchmark our competitive controller in a linearized Boeing 747 flight control system and show that it exhibits remarkable ``best-of-both-worlds" behavior, often beating standard $H_2$ and $H_{\infty}$ controllers on best-case and average-case input disturbances while maintaining a bounded loss in performance even in the worst-case. We also extend our competitive control framework to nonlinear systems using Model Predictive Control (MPC). Numerical experiments in a nonlinear system show that the competitive controller consistently outperforms standard $H_2$ and $H_{\infty}$ controllers across a wide variety of input disturbances, often by a large margin. This may be because the competitive controller, which is designed to adapt to arbitrary disturbance sequences, is better able to adapt to changing system dynamics; we plan to investigate this phenomenon more thoroughly in future work.

In this paper, we focus on designing online controllers which compete against clairvoyant offline controllers; it is natural to extend the idea of competitive control to other classes of comparator controllers. For example, it would be interesting to design distributed controllers which make decisions using only local information while competing against centralized controllers with a more global view. We anticipate that such distributed controllers could prove useful in a variety of networked control problems arising in congestion control, distributed resource allocation and smart grid.

\section{Appendix}

\subsection{Proof of Theorem \ref{competitive-controller-thm-fh}}

\begin{proof}
A state-space model for $F$ is given by $$\epsilon_{t+1} = A_t \epsilon_t + B_{u, t} u_t, \hspace{3mm} s_t = Q_t^{1/2} \epsilon_t.$$
Given this state-space model, we wish to obtain the factorization $\Delta \Delta^* = I + FF^*$ where $\Delta$ is causal. We interpret $I + FF^*$ as the covariance matrix of an appropriately defined random variable and use the Kalman filter to obtain a state-space model for $\Delta$. Suppose that $u$ and $v$ are
zero-mean random variables such that $\expect[uu^*] = I$, $\expect[vv^*] = I$ and $\expect[uv^*] = 0$.  Define $y = F u + v$; notice that $\expect[y y^*] = I + FF^*$. As is well-known in the signal processing community, the Kalman filter can be used to construct a causal matrix $\Delta$ such that $y = \Delta e$, where $e$ is a zero-mean random variable such that $\expect[e e^*] = I$; this is the so-called ``whitening" property of the Kalman filter. Notice that since $y = Fu + v$, $\expect[y y^*] = I + FF^*$; on the other hand, $y = \Delta e$, so $\expect[y y^*] = \Delta \Delta^*$. Therefore $I + FF^* = \Delta \Delta^*$, as desired. 

Using the Kalman filter as described in Theorem 9.2.1 in \cite{kailath2000linear}, we obtain a state-space model for $\Delta$:
\begin{equation} \label{Delta-ss}
\eta_{t+1} = A_t \eta_t + K_{t} \Sigma^{1/2}_{t} e_t, \hspace{3mm} y_t = Q_t^{1/2} \eta_t +  \Sigma_t^{1/2}e_t,
\end{equation}
where we define 
\begin{equation} \label{K-and-Sigma}
K_{t} = A_t P_t Q_t^{1/2} \Sigma_t^{-1}, \hspace{3mm} \Sigma_t = I + Q_t^{1/2}P_t Q_t^{1/2},
\end{equation}
and $P_t$ is defined recursively as $$P_{t+1} = A_tP_tA_t^* + B_{u, t}B_{u, t}^* - K_t \Sigma_tK_t^*$$
where we initialize $P_0 = 0$. 

Now that we have state-space models for $F$ and $\Delta$, we can form a state-space model for the overall system (\ref{transformed-dynamics}).
Letting $\alpha_t = \epsilon_t + \eta_t$, we see that a state-space model for this system is  $$\alpha_{t+1} = A_t \alpha_t + B_{u, t} u_t + K_t\Sigma_t^{1/2}w'_t, \hspace{3mm} s_t = Q_t^{\frac{1}{2}}\alpha_t + \Sigma_t^{\frac{1}{2}}w_t'.$$
This system can be rewritten as 
\begin{equation} \label{final-ss-2}
\xi_{t+1} = \hat{A}_t \xi_t + \hat{B}_{u, t} u_t + \hat{B}_{w, t} w_{t+1}', \hspace{3mm}
s_t = \hat{Q}_t^{1/2} \xi_t, 
\end{equation}
where we define 
$$\hat{A}_t = \begin{bmatrix} A_t & K_t \Sigma_t^{1/2} \\ 0 & 0 \end{bmatrix}, \hspace{3mm} \hat{B}_{u, t} = \begin{bmatrix} B_{u, t} \\ 0 \end{bmatrix}, \hspace{3mm} \hat{B}_{w, t} = \begin{bmatrix} 0 \\ I \end{bmatrix}, $$ $$ \hat{Q}_t = \begin{bmatrix} Q_t & Q_t^{1/2} \Sigma_t^{1/2} \\ \Sigma_t^{1/2}Q_t^{1/2} & \Sigma_t \end{bmatrix}
$$
and we initialize $\xi_0 = 0$. Recall that our goal is to find a controller in the synthetic system (\ref{final-ss-2}) such that $ALG(w) < \gamma^2 \|w'\|_2^2$ for all disturbances $w'$, or to determine whether no such controller exists; such a controller has competitive ratio at most $\gamma^2$ in the original system (\ref{evolution-eq}). Theorem \ref{hinf-suboptimal-controller-thm-fh} gives necessary and sufficient conditions for the existence of such a controller, along with an explicit state-space description of the controller, if it exists. 

We emphasize that the driving disturbance in the synthetic system (\ref{final-ss-2}) is not $w$, but rather the synthetic disturbance $w' = \Delta^{-1} G w$. Notice that $\Delta^{-1}G$ is strictly causal, since $\Delta^{-1}$ is causal and $G$ is strictly causal. Exchanging inputs and outputs in (\ref{Delta-ss}), we see that a state-space model for $\Delta^{-1}$ is $$\eta_{t+1} = (A_t - K_tQ_t^{1/2})\eta_t + K_ty_t, \hspace{3mm} e_t = \Sigma_t^{-1/2}(y_t - Q_t^{1/2} \eta_t).$$ A state-space model for $G$ is $$
\delta_{t+1} = A_t \delta_t + B_{w, t} w_t, \hspace{3mm} s_t = Q_t^{1/2} \delta_t.$$ 
Equating $s$ and $y$, we see that a state-space model for $\Delta^{-1}G$ is  
$$ \begin{bmatrix} \eta_{t+1} \\ \delta_{t+1} \end{bmatrix} = \begin{bmatrix} A_t - K_t Q_t^{1/2} & K_tQ_t^{1/2}  \\ 0 & A_t  \end{bmatrix} \begin{bmatrix} \eta_t \\ \delta_t \end{bmatrix} + \begin{bmatrix} 0 \\ B_{w, t}  \end{bmatrix} w_t, $$ $$ e_t =  \Sigma_t^{-1/2}Q_t^{1/2}( \delta_t - \eta_t). 
$$
Setting $\nu_t = \delta_t - \eta_t$ and simplifying, we see that a minimal representation for $w'$ is
$$\nu_{t+1} = (A_t - K_t Q_t^{1/2})\nu_t + B_{w, t}w_t, \hspace{3mm} w'_t = \Sigma_t^{-1/2} Q_t^{1/2} \nu_t.$$
We reiterate that $w'$ is a strictly casual function of $w$; in particular, $w_{t+1}'$ depends only on $w_0, w_1, \ldots, w_t$. 
\end{proof}

\subsection{Proof of Theorem \ref{competitive-controller-thm-ih}}

\begin{proof}
Taking the $z$-transform of the linear evolution equations $$x_{t+1} = A x_t + B_u u_t + B_w w_t, \hspace{3mm} s_t = Q^{1/2} x_t$$ we obtain $$zx(z) = Ax(z) + B_u u(z) +  B_w w(z), \hspace{3mm} s(z) = Q^{1/2} x(z)$$ Letting $F(z)$ and $G(z)$ be the transfer operators mapping $u(z)$ and $w(z)$ to $s(z)$, respectively, we see that $$F(z) = Q^{1/2}(zI - A)^{-1}B_u $$ and $$G(z) = Q^{1/2}(zI - A)^{-1}B_w.$$ Our goal is to obtain a canonical factorization $$I + F(z) F(z^{-*})^* = \Delta(z) \Delta(z^{-*})^*.$$ With this factorization, we can easily recover the optimal infinite-horizon competitive controller; it simply the $H_{\infty}$-optimal infinite-horizon controller in the system whose dynamics in the frequency domain are 
\begin{equation} \label{freq-domain-dynamics-delta}
s(z) = F(z)u(z) + \Delta(z)w'(z),
\end{equation}
where the synthetic disturbance $w'$ is 
\begin{equation} \label{freq-domain-disturbance}
w'(z) = \Delta^{-1}(z)G(z)w(z).
\end{equation}

Before we factor $ I + F(z) F(z^{-*})^*$, we state a key identity which plays a pivotal role in the factorization: for all Hermitian matrices $P$, we have
$$\begin{bmatrix} Q^{1/2}(zI - A)^{-1} & I \end{bmatrix} \Omega(P) \begin{bmatrix} (z^{-*}I - A)^{-*}Q^{1/2} \\ I \end{bmatrix} = 0, $$ where we define $$\Omega(P) = \begin{bmatrix} -P + APA^* & APQ^{1/2} \\ Q^{1/2}PA^* & Q^{1/2}PQ^{1/2} \end{bmatrix}. $$ This identity is easily verified via direct calculation.

We expand $ I + F(z) F(z^{-*})^*$ as
$$  \begin{bmatrix} Q^{1/2}(zI - A)^{-1} & I \end{bmatrix} \begin{bmatrix} B_uB_u^* & 0 \\ 0 & I \end{bmatrix} \begin{bmatrix} (z^{-*}I - A)^{-*}Q^{1/2} \\ I \end{bmatrix}.$$
Applying the identity, we see that this equals
$$ \begin{bmatrix} Q^{1/2}(zI - A)^{-1} & I \end{bmatrix} \Lambda(P) \begin{bmatrix} (z^{-*}I - A)^{-*}Q^{1/2} \\ I \end{bmatrix},$$
where $P$ is an arbitrary Hermitian matrix and we define $$\Lambda(P) = \begin{bmatrix} B_uB_u^* - P + APA^* & APQ^{1/2} \\ Q^{1/2}PA^* & I + Q^{1/2}PQ^{1/2} \end{bmatrix}.$$ Notice that the $\Lambda(P)$ can be factored as $$\begin{bmatrix} I & K(P) \\ 0 & I \end{bmatrix} \begin{bmatrix} \Gamma(P) & 0 \\ 0 & \Sigma(P) \end{bmatrix} \begin{bmatrix} I & 0 \\ K^*(P) & I \end{bmatrix}, $$
where we define $$\Gamma(P) = B_u B_u^* - P + APA^* - K\Sigma K^*(P),$$ 
\begin{equation} \label{sigma-k-ih-definition}
K(P) = APQ^{1/2}\Sigma^{-1}, \hspace{3mm} \Sigma(P) = I + Q^{1/2}PQ^{1/2}.
\end{equation}
By assumption, $(A, B_u)$ is stabilizable and $(A, Q^{1/2})$ is detectable, therefore the Ricatti equation $\Gamma(P) = 0$ has a unique stabilizing solution (see, e.g. Theorem E.6.2 in \cite{kailath2000linear}).  Suppose $P$ is chosen to be this solution, and define $K = K(P)$, $\Sigma = \Sigma(P)$. We immediately obtain the canonical factorization $$ I + F(z) F(z^{-*})^* = \Delta(z) \Delta^*(z^{-*}), $$ where we define 
\begin{equation} \label{delta-definition}
\Delta(z) = (I + Q^{1/2}(zI - A)^{-1}K )\Sigma^{1/2}.
\end{equation}
Define $$\hat{A} = \begin{bmatrix} A & K\Sigma^{1/2} \\ 0 & 0 \end{bmatrix}, \hspace{3mm} \hat{B}_u = \begin{bmatrix} B_u \\ 0 \end{bmatrix}, \hspace{3mm} \hat{B}_w = \begin{bmatrix} 0 \\ I \end{bmatrix}. $$ Notice that $\hat{\Delta}(z) = z^{-1}\Delta(z)$ can be cleanly expressed as 
\begin{equation} \label{hat-delta-z-transform}
\hat{\Delta}(z) = \begin{bmatrix} Q^{1/2} & \Sigma^{1/2} \end{bmatrix} (zI - \hat{A} )^{-1}\hat{B}_w.
\end{equation}
Similarly, $F(z)$ can be written as 
\begin{equation} \label{F-z-transform}
F(z) = \begin{bmatrix} Q^{1/2} & \Sigma^{1/2} \end{bmatrix}(zI - \hat{A} )^{-1}\hat{B}_u.
\end{equation}
We can rewrite the frequency domain dynamics (\ref{freq-domain-dynamics-delta}) in terms of $\hat{\Delta}(z)$:
\begin{equation} \label{freq-domain-dynamics-hat-delta}
s(z) = F(z)u(z) + \hat{\Delta}(z)(zw'(z)).
\end{equation}

It is easy to check that the stabilizability of $(A, B_u)$ implies the stabilizability of $(\hat{A}, \hat{B_u})$. Similarly, the $(A, Q^{1/2})$ is detectable and hence unit circle observable, which implies that $(\hat{A}, \begin{bmatrix} Q^{1/2} & \Sigma^{1/2} \end{bmatrix})$ is also unit circle observable. Applying Theorem \ref{hinf-ih-controller-thm} to the system (\ref{freq-domain-dynamics-hat-delta}) using the models for $\hat{\Delta}(z)$ and $F(z)$ given in (\ref{hat-delta-z-transform}) and (\ref{F-z-transform}), respectively, we obtain necessary and sufficient conditions for the existence of a competitive controller at level $\gamma$. We see that a frequency domain model $K(z)$ of the competitive controller at level $\gamma$ (if one exists) is given by $$-\hat{H}^{-1}\hat{B}_u^*\hat{P}\left[I + \hat{A}(zI - \hat{A}_2)^{-1}(I - \hat{B}_u^*\hat{H}^{-1}\hat{B}_u^*\hat{P}) \right]\hat{B}_w,$$ where we define $ \hat{H} = I + \hat{B}_u^*\hat{P}\hat{B}_u$, $\hat{A}_2 = \hat{A} - \hat{B}_u \hat{H}^{-1}\hat{B}_u^*\hat{P}\hat{A},$  and $\hat{P}$ is the solution of the Ricatti equation $$\hat{P} = \hat{Q} + \hat{A}^*\hat{P}\hat{A} -  \hat{A}^*\hat{P}\tilde{B} \tilde{H}^{-1} \tilde{B}^* \hat{P}\hat{A},$$ where we define $$\hat{Q} = \begin{bmatrix} Q^{1/2} \\ \Sigma^{1/2} \end{bmatrix} \begin{bmatrix} Q^{1/2} & \Sigma^{1/2} \end{bmatrix}, $$ $$ \tilde{B} = \begin{bmatrix} \hat{B}_{u} & \hat{B}_{w} \end{bmatrix}, \hspace{3mm} \tilde{H} = \begin{bmatrix} I & 0 \\ 0 & -\gamma^2 I \end{bmatrix}  + \tilde{B}^* \hat{P} \tilde{B}.$$

Translating this result back into time domain, we obtain a state-space model of the infinite-horizon controller: 
\begin{eqnarray*}
\xi_{t+1} &=& \hat{A} \xi_t + \hat{B}_u u_t + \hat{B}_w w'_{t+1} \\
u_t &=& -\hat{H}^{-1}\hat{B}_u^*\hat{P}(\hat{A}\xi_t + \hat{B}_w w'_{t+1}).
\end{eqnarray*}
Note that this system is driven by $w'_{t+1}$, not $w_t$, since the driving disturbance in (\ref{freq-domain-dynamics-hat-delta}) is $zw'(z)$.

We now construct the synthetic disturbance $w'$. Recall that $w'(z) = \Delta^{-1}(z)G(z)w(z)$. We have $$ \Delta^{-1}(z) = \Sigma^{-1/2} \left(I - Q^{1/2}(zI - (A - KQ^{1/2}))^{-1}K \right), $$ $$G(z) = Q^{1/2}(zI - A)^{-1}B_w.$$ We note that $A - KQ^{1/2}$ is stable and hence $\Delta^{-1}(z)$ is causal and bounded since its poles are strictly contained in the unit circle.
A state-space model for $w'$ is
$$ \begin{bmatrix} \eta_{t+1} \\ \delta_{t+1} \end{bmatrix} = \begin{bmatrix} A - K Q^{1/2} & KQ^{1/2}  \\ 0 & A  \end{bmatrix} \begin{bmatrix} \eta_t \\ \delta_t \end{bmatrix} + \begin{bmatrix} 0 \\ B_{w}  \end{bmatrix} w_t, $$ $$ w_t' =  \Sigma^{-1/2}Q^{1/2}( \delta_t - \eta_t). $$
Setting $\nu_t = \delta_t - \eta_t$ and simplifying, we see that a minimal representation for $w'$ is given by
$$\nu_{t+1} = (A - K Q^{1/2})\nu_t + B_{w}w_t, \hspace{3mm} w' = \Sigma^{-1/2} Q^{1/2} \nu.$$
We reiterate that $w'$ is a strictly casual function of $w$; in particular, $w_{t+1}'$ depends only on $w_0, w_1, \ldots, w_t$.
\end{proof}

%\newpage
\bibliography{main}

\begin{thebibliography}{10}

\bibitem{agarwal2019online}
N.~Agarwal, B.~Bullins, E.~Hazan, S.~M. Kakade, and K.~Singh.
\newblock Online control with adversarial disturbances.
\newblock {\em arXiv preprint arXiv:1902.08721}, 2019.

\bibitem{borodin2005online}
A.~Borodin and R.~El-Yaniv.
\newblock {\em Online computation and competitive analysis}.
\newblock cambridge university press, 2005.

\bibitem{chen2018smoothed}
N.~Chen, G.~Goel, and A.~Wierman.
\newblock Smoothed online convex optimization in high dimensions via online
  balanced descent.
\newblock In {\em Conference On Learning Theory}, pages 1574--1594. PMLR, 2018.

\bibitem{doyle1978guaranteed}
J.~C. Doyle.
\newblock Guaranteed margins for lqg regulators.
\newblock {\em IEEE Transactions on automatic Control}, 23(4):756--757, 1978.

\bibitem{foster2020logarithmic}
D.~J. Foster and M.~Simchowitz.
\newblock Logarithmic regret for adversarial online control.
\newblock {\em arXiv preprint arXiv:2003.00189}, 2020.

\bibitem{goel2017thinking}
G.~Goel, N.~Chen, and A.~Wierman.
\newblock Thinking fast and slow: Optimization decomposition across timescales.
\newblock In {\em 2017 IEEE 56th Annual Conference on Decision and Control
  (CDC)}, pages 1291--1298. IEEE, 2017.

\bibitem{goel2020power}
G.~Goel and B.~Hassibi.
\newblock The power of linear controllers in lqr control.
\newblock {\em arXiv preprint arXiv:2002.02574}, 2020.

\bibitem{goel2021regret}
G.~Goel and B.~Hassibi.
\newblock Regret-optimal estimation and control.
\newblock {\em arXiv preprint arXiv:2106.12097}, 2021.

\bibitem{goel2021measurement}
G.~Goel and B.~Hassibi.
\newblock Regret-optimal measurement-feedback control.
\newblock In {\em Learning for Dynamics and Control}, pages 1270--1280. PMLR,
  2021.

\bibitem{goel2019beyond}
G.~Goel, Y.~Lin, H.~Sun, and A.~Wierman.
\newblock Beyond online balanced descent: An optimal algorithm for smoothed
  online optimization.
\newblock In {\em Advances in Neural Information Processing Systems}, pages
  1875--1885, 2019.

\bibitem{goel2019online}
G.~Goel and A.~Wierman.
\newblock An online algorithm for smoothed regression and lqr control.
\newblock {\em Proceedings of Machine Learning Research}, 89:2504--2513, 2019.

\bibitem{gradu2020adaptive}
P.~Gradu, E.~Hazan, and E.~Minasyan.
\newblock Adaptive regret for control of time-varying dynamics.
\newblock {\em arXiv preprint arXiv:2007.04393}, 2020.

\bibitem{hassibi1999indefinite}
B.~Hassibi, A.~H. Sayed, and T.~Kailath.
\newblock {\em Indefinite-quadratic estimation and control: a unified approach
  to H 2 and H-infinity theories}.
\newblock SIAM, 1999.

\bibitem{hazan2020nonstochastic}
E.~Hazan, S.~Kakade, and K.~Singh.
\newblock The nonstochastic control problem.
\newblock In {\em Algorithmic Learning Theory}, pages 408--421. PMLR, 2020.

\bibitem{boeing747}
J.~Hong, N.~Moehle, and S.~Boyd.
\newblock Lecture notes in ”introduction to matrix methods", 2021.

\bibitem{kailath2000linear}
T.~Kailath, A.~H. Sayed, and B.~Hassibi.
\newblock {\em Linear estimation}.
\newblock Prentice Hall, 2000.

\bibitem{sabag2021regret}
O.~Sabag, G.~Goel, S.~Lale, and B.~Hassibi.
\newblock Regret-optimal full-information control.
\newblock {\em arXiv preprint arXiv:2105.01244}, 2021.

\bibitem{shi2020online}
G.~Shi, Y.~Lin, S.-J. Chung, Y.~Yue, and A.~Wierman.
\newblock Online optimization with memory and competitive control.
\newblock {\em arXiv e-prints}, pages arXiv--2002, 2020.

\bibitem{zhao2021non}
P.~Zhao, Y.-X. Wang, and Z.-H. Zhou.
\newblock Non-stationary online learning with memory and non-stochastic
  control.
\newblock {\em arXiv preprint arXiv:2102.03758}, 2021.

\end{thebibliography}

\end{document}